\documentclass[journal]{IEEEtran}

\usepackage{amsmath,amssymb,amsfonts}
\usepackage{graphicx}
\usepackage{cite}
\usepackage{subfigure}
\usepackage{algorithm}
\usepackage{algorithmic}
\usepackage{multirow}
\usepackage{booktabs}
\usepackage{tikz}
\usetikzlibrary{arrows.meta, calc, shapes.geometric, positioning, fit}
\usepackage{pifont}
\usepackage[export]{adjustbox}
\usepackage{tabularray}
\usepackage{colortbl}
\usepackage{array}
\usepackage{cuted}
\usepackage{capt-of}
\usepackage[colorlinks,linkcolor=blue,anchorcolor=blue,citecolor=blue]{hyperref}

\definecolor{tableTitleColor}{rgb}{0.8,0.8,0.8}
\newcommand{\modelcodeicon}[2]{%
  \begin{scope}[shift={#1}, scale=#2]
    \draw[black!72, line width=0.78pt, rounded corners=1pt]
      (-0.14,-0.46) rectangle (0.42,0.50);
    \draw[black!62, line width=0.62pt, line cap=round]
      (-0.02,0.28) -- (0.30,0.28)
      (-0.02,0.14) -- (0.30,0.14)
      (-0.02,0.00) -- (0.22,0.00);
    \coordinate (gearc) at (-0.22,-0.28);
    \foreach \a in {0,45,...,315} {
      \draw[black!72, line width=1.00pt, line cap=round]
        ($(gearc)+(\a:0.20)$) -- ($(gearc)+(\a:0.29)$);
    }
    \draw[black!72, fill=white, line width=0.78pt] (gearc) circle (0.22);
    \draw[black!72, fill=white, line width=0.72pt] (gearc) circle (0.09);
  \end{scope}%
}

\newcommand{\taskdescicon}[2]{%
  \begin{scope}[shift={#1}, scale=#2]
    \draw[black!65, line width=0.72pt, rounded corners=1pt]
      (-0.42,-0.50) -- (-0.42,0.50) -- (0.34,0.50) --
      (0.34,-0.50) -- (-0.24,-0.50) -- cycle;
    \draw[black!65, line width=0.72pt, line join=round]
      (-0.42,-0.26) -- (-0.24,-0.26) -- (-0.24,-0.50);
    \draw[black!60, line width=0.72pt, line cap=round]
      (-0.11,0.17) -- (0.10,0.30)
      (-0.11,0.17) -- (0.10,0.02);
    \foreach \p in {(-0.20,0.17),(0.16,0.34),(0.16,-0.01)} {
      \draw[black!70, fill=gray!10, line width=0.72pt] \p circle (0.075);
    }
    \draw[black!52, line width=0.64pt, line cap=round]
      (-0.02,-0.22) -- (0.22,-0.22)
      (-0.10,-0.34) -- (0.22,-0.34);
  \end{scope}%
}

\newcommand{\solutioncheckicon}[2]{%
  \begin{scope}[shift={#1}, scale=#2]
    \draw[black!65, line width=0.72pt, rounded corners=1pt]
      (-0.38,-0.50) rectangle (0.22,0.50);
    \foreach \y in {0.26,0,-0.26} {
      \draw[black!65, line width=0.60pt, rounded corners=0.5pt]
        (-0.28,\y-0.055) rectangle (-0.17,\y+0.055);
      \draw[black!52, line width=0.60pt, line cap=round]
        (-0.10,\y+0.025) -- (0.12,\y+0.025)
        (-0.10,\y-0.045) -- (0.05,\y-0.045);
    }
    \draw[black!65, fill=gray!10, line width=0.72pt]
      (0.30,-0.02) circle (0.27);
    \draw[black!78, line width=0.86pt, line cap=round, line join=round]
      (0.17,-0.02) -- (0.27,-0.13) -- (0.46,0.11);
    \draw[black!65, line width=0.72pt, rounded corners=1pt]
      (0.24,-0.62) rectangle (0.36,-0.32);
  \end{scope}%
}

\newcommand{\controllericon}[2]{%
  \begin{scope}[shift={#1}, scale=#2]
    \path[draw=violet!70!black, fill=violet!8, line width=0.82pt]
      (-0.46,-0.10)
        .. controls (-0.42,0.15) and (-0.27,0.26) .. (-0.10,0.18)
        -- (0.10,0.18)
        .. controls (0.27,0.26) and (0.42,0.15) .. (0.46,-0.10)
        .. controls (0.54,-0.34) and (0.31,-0.45) .. (0.16,-0.25)
        -- (-0.16,-0.25)
        .. controls (-0.31,-0.45) and (-0.54,-0.34) .. (-0.46,-0.10)
        -- cycle;
    \draw[violet!70!black, line width=0.78pt, line cap=round]
      (-0.30,-0.02) -- (-0.14,-0.02)
      (-0.22,-0.10) -- (-0.22,0.06);
    \draw[violet!70!black, fill=white, line width=0.70pt]
      (0.21,0.02) circle [radius=0.045]
      (0.33,-0.07) circle [radius=0.045];
    \draw[violet!58!black, fill=white, line width=0.64pt]
      (-0.04,-0.04) circle [radius=0.035]
      (0.07,-0.04) circle [radius=0.035];
  \end{scope}%
}

\begin{document}

\title{Vehicle Routing Problem Meets Large Language Models: An Overview and Perspectives}

\author{Xianchao Xiu,~\IEEEmembership{Member,~IEEE}, Chong Shen, Yanjiao Zhu, and Wanquan Liu,~\IEEEmembership{Senior Member,~IEEE}%
\thanks{This work was supported in part by the National Natural Science Foundation of China under Grants 12371306. (\textit{Corresponding author: Wanquan Liu}.)}
\thanks{Xianchao Xiu is with the School of Mechatronic Engineering and Automation, Shanghai University, Shanghai 200444, China (e-mail: xcxiu@shu.edu.cn).}%
\thanks{Chong Shen is with the Sino-European School of Technology, Shanghai University, Shanghai 200444, China (e-mail: liuwq63@mail.sysu.edu.cn).}%
\thanks{Yanjiao Zhu and Wanquan Liu are with the School of Intelligent Systems Engineering, Sun Yat-sen University, Shenzhen 518107, China (e-mail: liuwq63@mail.sysu.edu.cn).}%
}

\maketitle

\begin{abstract}
The vehicle routing problem (VRP) is a central optimization problem in artificial intelligence, logistics automation, 
transportation scheduling, and industrial decision-making. VRP and its variants are NP-hard, and practical routing tasks often combine time windows, 
vehicle capacities, pickup-and-delivery relations, dynamic requests, and other operational constraints, making both modeling and solving difficult. 
Large language models (LLMs) provide a flexible interface for routing optimization by processing natural-language requirements, generating code, reasoning over constraints, 
and interacting with external tools. This survey reviews LLM-driven research on VRP, covering the basic definition, main variants, major solver families, 
and LLM concepts needed for this topic. Existing studies are organized into three roles: modelers translate natural-language requirements into constraints and modeling code; 
designers generate heuristics, operators, or route plans; and coordinators organize tool calls, multi-agent collaboration, and connections with neural solvers. 
The survey also reviews standard benchmarks, real or near-real operational datasets, LLM-oriented evaluation frameworks, and two comparative experiments. 
The goal is to clarify current progress in LLM-assisted routing optimization and provide a structured reference for intelligent decision-making, 
advanced manufacturing, and industrial automation.
\end{abstract}

\begin{IEEEkeywords}
Vehicle routing problem, large language models, combinatorial optimization, heuristic design, multi-agent systems.
\end{IEEEkeywords}

\section{Introduction}
\label{sec:introduction}

\IEEEPARstart{T}{HE} vehicle routing problem (VRP) describes routing and resource-allocation decisions in logistics distribution, 
instant delivery, warehouse replenishment, ride-hailing, and vehicle dispatching. A routing system must decide which vehicles serve which customers, 
the order of visits, and how service deadlines and resource limits are satisfied. These decisions directly affect transportation cost,
 fulfillment efficiency, and system reliability. In practice, VRP instances are rarely governed by a single constraint. Vehicle capacity, 
 time windows, multiple depots, charging or refueling, pickup-and-delivery pairing, stochastic demand, and dynamic disruptions may appear together. 
 As a result, VRP is not only a difficult optimization problem but also a demanding modeling and deployment problem. 
 New combinations of constraints often require fresh formulations, custom operators, parameter tuning, and careful solver adaptation, 
 which limits the transferability of existing routing systems.

VRP solvers have developed along three broad lines. Exact methods use mathematical programming, branch-and-price, cutting planes, 
dynamic programming, and related techniques to prove optimality or provide rigorous bounds. Their computational cost, however, 
can grow sharply with instance size and constraint complexity. Heuristic and metaheuristic methods, including constructive rules, 
local search, destroy-and-repair schemes, and population-based evolution, 
trade optimality guarantees for high-quality feasible solutions within practical time budgets. 
NCO and other learning-enhanced methods use data distributions, policy learning, and neural representations to improve route construction or 
guide search \cite{bengio2021machine}, but they still struggle with feasibility guarantees, scale transfer, and unfamiliar constraint combinations.

\begin{figure*}[!t]
\centering
\includegraphics[width=0.88\textwidth]{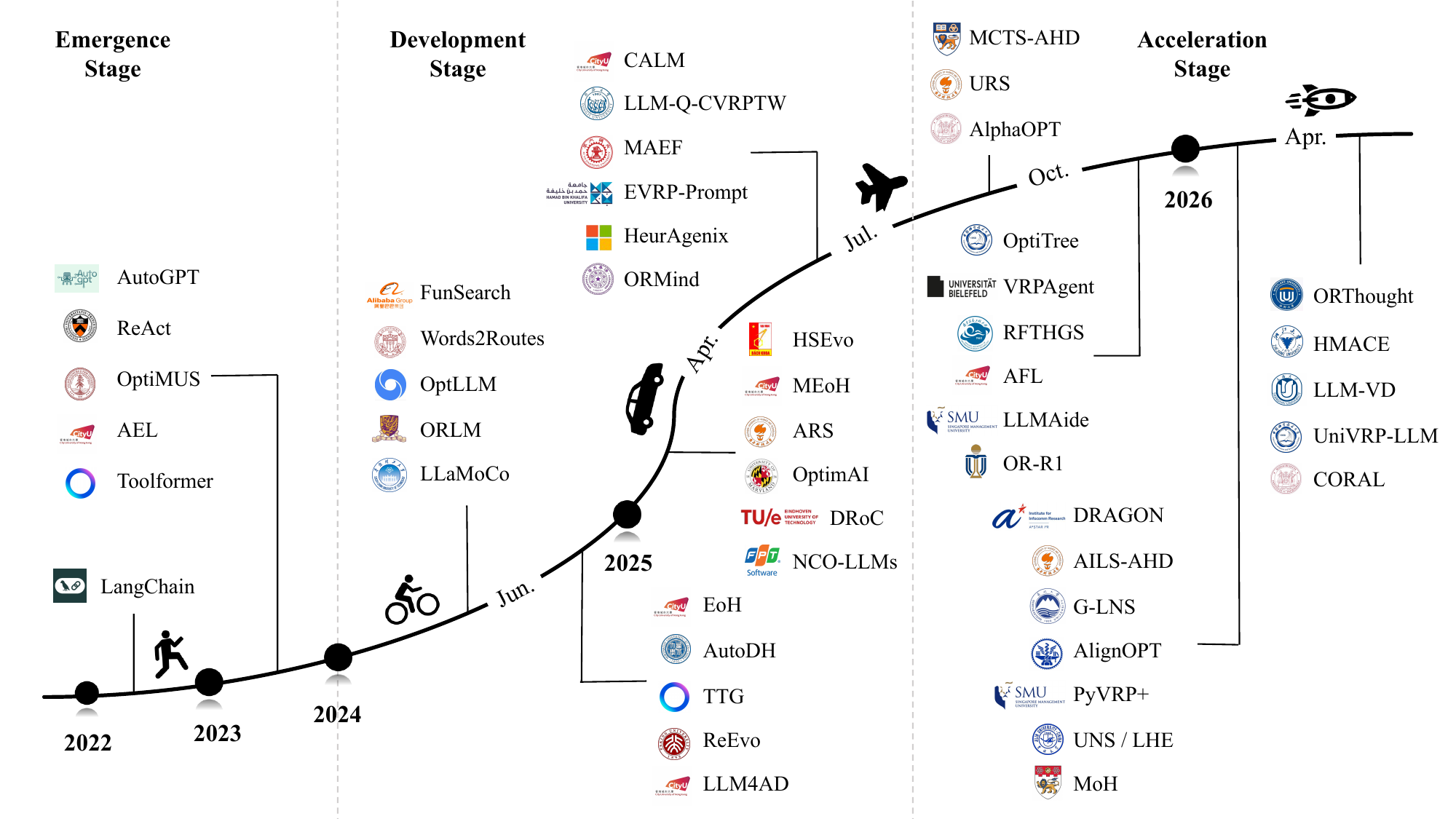}
\caption{Timeline of LLM-driven VRP research.}
\label{fig:vrp-solver-timeline}
\end{figure*}

Large language models (LLMs) have shown strong capabilities in natural-language understanding \cite{brown2020language}, code generation \cite{chen2021evaluating}, 
and algorithm design \cite{liu2024llm_survey_algorithm_design}. These capabilities create a flexible interface between business descriptions 
and optimization pipelines. In LLM-driven VRP research, models are used to parse natural-language constraints, generate mathematical-modeling code, 
propose heuristics or neighborhood operators, call external tools, and coordinate multi-agent workflows. As shown in Fig.~\ref{fig:vrp-solver-timeline}, 
early studies were closely tied to general agent and tool-use frameworks. Later work moved toward automatic modeling, heuristic design, 
and direct route reasoning, and recent research has increasingly examined multi-agent collaboration and connections with neural solvers.

Recent surveys have examined LLMs for combinatorial optimization \cite{xiu2026llmor_survey}, 
operations research and optimization \cite{wang2025llmor_survey}, metaheuristic optimization \cite{ghanbarzadeh2025structured}, 
and multimodal foundation models \cite{albalkhi2026routeoptimization}. However, existing surveys rarely analyze the different roles of LLMs in modeling, 
solution generation, and workflow coordination from the perspective of the VRP solving process. 
VRP-related datasets and experimental results in related work also remain scattered.

This paper presents a systematic survey of LLM-driven VRP research. After reviewing more than 60 closely related studies, 
it groups existing work into three roles, modeler, designer, and coordinator, corresponding to problem modeling, solution generation, 
and solver workflow coordination. The survey also summarizes standard VRP benchmarks, operational-scenario datasets, and LLM evaluation frameworks, 
and discusses their applicability to different evaluation tasks. Finally, two comparative experiments compare different methods 
in terms of solution quality, solving-time cost, and scalability.

The remainder of this paper is organized as follows. Section~\ref{sec:background} introduces the basic concepts of VRP, 
its main variants, classical and learning-based solvers, and relevant background on LLMs. 
Sections~\ref{sec:modelers}--\ref{sec:coordinators} discuss LLMs from the roles of modeler, designer, and coordinator, 
covering automatic modeling and code generation, automated heuristic design and end-to-end route reasoning, tool-chain orchestration, 
multi-agent collaboration, and neural-solver integration. Section~\ref{sec:datasets} introduces VRP-related datasets. 
Section~\ref{sec:experiments} presents two comparative experiments. 
Section~\ref{sec:conclusion} concludes the paper and discusses future research directions.

\section{Background}
\label{sec:background}

Depending on constraints, objectives, and application scenarios, VRP has developed into a family of routing problems and a broad set of solution methods. 
This section reviews the basic formulation, common variants, main solver families, and LLM concepts needed to discuss modeling, 
heuristic generation, and workflow coordination.

\subsection{Vehicle Routing Problem}
\label{subsec:vrp}

\subsubsection{Classical Formulation}

VRP was first introduced by Dantzig and Ramser in 1959 and has since become a central problem in distribution optimization \cite{dantzig1959truck}. 
It can be regarded as an extension of the traveling salesman problem (TSP) from a single tour to vehicle scheduling and distribution. 
While TSP determines one visiting sequence, VRP must allocate customers among vehicles and design a route for each vehicle under capacity, 
service-time, and demand constraints. The basic task is to serve all customers from the available depots and fleet while minimizing total travel cost.

The capacitated VRP (CVRP) is the most basic and widely studied VRP variant. A fleet departs from a depot, 
serves geographically distributed customers, and returns to the depot without exceeding vehicle capacity \cite{letchford2015mcf_cvrp}. 
Let $V=\{0\}\cup V_c$ be the node set, where node $0$ is the depot and $V_c$ is the customer set. Let $A$ be the directed arc set without self-loops, 
and let $(i,j)\in A$ denote travel from node $i$ to node $j$. The arc cost is $c_{ij}$, the total demand of a customer subset $S$ is $q(S)$, 
and the vehicle capacity is $Q$. A compact CVRP formulation is

\begin{equation}
  \begin{aligned}
    \min_{x_{ij}} \quad & \sum_{(i,j)\in A} c_{ij} x_{ij} \\
    \mathrm{s.t.} \quad & x(\delta^{+}(i)) = 1, && \forall i \in V_c \\
    & x(\delta^{-}(i)) = 1, && \forall i \in V_c \\
    & x(\delta^{+}(S)) \geq \left\lceil \frac{q(S)}{Q} \right\rceil,
      && \forall S \subseteq V_c \\
    & x_{ij} \in \{0,1\}, && \forall (i,j)\in A
  \end{aligned}
  \label{eq:cvrp}
\end{equation}
Here, $x_{ij}$ indicates whether arc $(i,j)$ is used. The sets $\delta^{+}(i)$ and $\delta^{-}(i)$ are the outgoing and incoming arcs of node $i$. The degree constraints require each customer to be entered and left exactly once, while the subset constraints require at least $\left\lceil q(S)/Q \right\rceil$ routes to leave any customer subset $S$. The number of feasible solutions grows super-exponentially with the number of customers, and capacity restrictions further constrain both customer grouping and visiting order. Thus exact CVRP solving is NP-hard \cite{laporte2009fifty}, which is why CVRP has long been a core benchmark for routing optimization.

\subsubsection{Representative Variants}

Many VRP variants are obtained by adding operational constraints to the classical setting. The VRP with time windows (VRPTW) introduces service intervals, so visiting order and service schedules must be optimized together. The split delivery VRP (SDVRP) allows the demand of a customer to be served by several vehicles, relaxing the single-visit assumption of CVRP. The electric VRP (EVRP), multi-depot VRP (MDVRP), and VRP with backhauls (VRPB) emphasize charging or refueling, depot selection, and pickup operations on return legs, respectively. The open VRP (OVRP) removes the requirement that vehicles return to the depot after service, while the VRP with route length limit (VRPL) constrains the maximum distance or duration of each route to model driving range, working time, or service-radius restrictions. The periodic VRP (PVRP) extends routing decisions across multiple days or service periods, requiring both visit schedules and within-period routes to be determined.

These variants provide different task scenarios for LLM-assisted VRP solving. For example, the time-window constraints in VRPTW can be used to examine whether an LLM can reason about arrival time, waiting time, and service windows.

\subsection{VRP Solvers}
\label{subsec:vrp-solvers}

VRP solvers historically developed around exact optimization and heuristic search. Machine learning has recently been introduced for route construction, route improvement, and solver assistance. Exact methods seek optimal solutions or rigorous bounds and are most useful for small or structurally clear instances. Heuristics prioritize good feasible solutions under limited computation and are widely used in medium- and large-scale applications. Learning-based methods infer construction policies, improvement operators, or auxiliary decisions from historical instances, expert solutions, or search feedback, and usually complement rather than replace exact or heuristic frameworks.

\subsubsection{Exact Methods}

Exact methods first formulate VRP as an integer linear program and then reduce the feasible search space through branching, cutting planes, column generation, or their combinations. Branch-and-bound partitions the problem into subproblems and prunes them using upper and lower bounds. Branch-and-cut adds valid inequalities at branch nodes to remove fractional extreme points in linear relaxations and tighten the feasible region. Branch-and-price embeds column generation into branch search, introducing path or resource variables on demand to reduce the growth of the column space. Branch-price-and-cut further combines column generation and cutting planes in route- or task-set formulations, strengthening both bounds and variable spaces. Such methods can provide optimal solutions or rigorous bounds, but their computational cost usually increases rapidly with problem scale. As a result, they are mainly used for offline exact optimization on small and medium-sized instances.

\subsubsection{Heuristic and Metaheuristic Methods}

Heuristic methods aim to find high-quality feasible solutions under limited computational budgets. 
They can be broadly divided into constructive heuristics, improvement heuristics, and metaheuristics. 
Constructive heuristics start from an empty solution and gradually build complete routes. Representative examples include the savings algorithm, 
sweep algorithm, and insertion methods, which are often used to produce initial feasible solutions quickly. 
Improvement heuristics start from an existing solution and repeatedly improve route structures through neighborhood operations such as exchange, 
relocation, 2-opt, and 3-opt. Metaheuristics organize the search process at a higher level, with typical examples including tabu search, 
simulated annealing, genetic algorithms, and ant colony optimization \cite{laporte2009fifty}.

Among modern heuristic frameworks, large neighborhood search (LNS) flexibly combines destroy and repair operators and can be 
adapted to many VRP variants \cite{pisinger2007general}. Adaptive large neighborhood search (ALNS) further introduces adaptive operator selection 
and weight updates, making it more stable for large-scale and heavily constrained problems \cite{ropke2006adaptive}. 
Hybrid genetic search (HGS) combines genetic search, local search, and diversity management, 
and has become a strong baseline across multiple VRP benchmarks \cite{vidal2013hybrid}. However, 
the performance of these methods still depends heavily on manual operator design, search strategy construction, 
and parameter tuning. When new constraint combinations or problem variants appear, substantial transfer and redesign costs are often required. 
This limitation has motivated learning-based methods that seek to reduce manual engineering effort.

\subsubsection{Machine-Learning-Based Methods}
\label{subsec:ml-solvers}

The core idea of machine-learning-based routing methods is to use neural networks to learn construction, improvement, 
or assistance strategies from data, expert solutions, or interactive feedback. NCO is a major line of recent work. 
Since VRP instances can naturally be represented as graphs containing depots, customers, and edge weights, graph neural networks, 
graph attention networks, and Transformer encoders are commonly used to learn node, edge, and global instance representations. 
Training can be conducted through supervised learning by imitating optimal or expert solutions, 
or through reinforcement learning by directly using route cost and constraint penalties as feedback. 
Policy gradient and actor-critic methods are common reinforcement learning choices.

Several end-to-end neural solvers have become important milestones. 
Pointer Networks \cite{vinyals2015pointer} introduced attention mechanisms for variable-length sequence decisions and provided an early framework 
for directly outputting visit orders. The Attention Model (AM) \cite{kool2018attention} introduced multi-head attention and graph embeddings, 
improving the representation of customers, depots, and edge relations. POMO \cite{kwon2020pomo} exploits multiple equivalent optimal solutions 
within the same instance to enhance policy-gradient training and improve sample efficiency and solution quality. Sym-NCO \cite{kim2022symnco} uses 
symmetries such as rotation and reflection to construct regularization mechanisms, 
improving the generalization of deep-reinforcement-learning-based NCO without relying on problem-specific expert knowledge. 
LEHD \cite{luo2023lehd} improves generalization on larger TSP and CVRP instances through a light encoder and heavy decoder architecture. 
RouteFinder \cite{berto2025routefinder} further constructs a unified neural routing solver for multiple VRP variants, 
aiming to improve shared representation and transfer across routing tasks.

Beyond directly generating routes, machine learning can also serve as an enhancement module inside traditional solvers. 
In exact methods, learning models can assist branching-variable selection, node selection, cut selection, or column-generation pricing. 
In heuristic and metaheuristic methods, they can be used for initial-solution generation, neighborhood selection, destroy-repair operator selection, 
and parameter adjustment. Prior work has used graph convolutional networks to learn branching decisions in mixed-integer programming 
(MIP) \cite{gasse2019gcnn}, and neural networks have also been embedded into the repair process of LNS \cite{hottung2020nlNS}. 
Machine learning can improve efficiency on repeated solving and in-distribution instances, but its generalization across scales,
distributions, and complex constraints remains limited. This limitation creates demand for more flexible optimization interfaces that can 
handle natural-language constraints, code generation, and external tool use.

\subsection{Large Language Models}
\label{subsec:llms}

\subsubsection{Basic Principles}

LLMs are primarily built on the Transformer architecture \cite{zhou2020progress}. Its multi-head self-attention mechanism captures dependencies among positions in a sequence \cite{vaswani2017attention}, making Transformers effective for long-context processing. Given an input representation matrix $\mathbf{X}$, self-attention first projects it into query, key, and value matrices:

\begin{equation}
  \mathbf{Q}=\mathbf{X}\mathbf{W}_Q,\quad \mathbf{K}=\mathbf{X}\mathbf{W}_K,\quad \mathbf{V}=\mathbf{X}\mathbf{W}_V
  \label{eq:qkv_projection}
\end{equation}
where $\mathbf{W}_Q$, $\mathbf{W}_K$, and $\mathbf{W}_V$ are parameter matrices. Attention weights are then computed from similarities between queries and keys and used to aggregate the value matrix:

\begin{equation}
  \operatorname{Attention}(\mathbf{Q},\mathbf{K},\mathbf{V})
  =\operatorname{softmax}\left(\frac{\mathbf{Q}\mathbf{K}^{\top}}{\sqrt{d_k}}\right)\mathbf{V}
  \label{eq:self_attention}
\end{equation}
Here, $d_k$ is the dimension of the key vectors. The product $\mathbf{Q}\mathbf{K}^{\top}$ measures correlations among position representations, and the normalized attention weights aggregate $\mathbf{V}$ to form context-aware representations. Multi-head attention performs this computation in several representation subspaces in parallel, which improves the ability of the model to capture diverse semantic relations and long-range dependencies.

Most modern LLMs use decoder-only autoregressive architectures. During pretraining, they learn to predict the next token from large text corpora and thereby acquire language patterns, knowledge associations, and contextual dependencies \cite{min2023recent}. Supervised fine-tuning (SFT) \cite{dong2023sft} and preference alignment \cite{ouyang2022training} then improve instruction following, reasoning, and code generation. Representative models include the GPT series of OpenAI \cite{openai2023gpt4}, the Qwen series of Alibaba \cite{qwen2024qwen25}, and the DeepSeek series \cite{deepseek2025r1}. For resource-limited or domain-specific settings, parameter-efficient fine-tuning (PEFT) is important. Low-rank adaptation (LoRA) \cite{hu2022lora}, for example, freezes pretrained weights and adds low-rank trainable updates to selected linear layers, enabling task adaptation with relatively few new parameters. This makes domain adaptation to VRP-related data more practical.

\subsubsection{Functional Roles in VRP}
\label{subsec:llm-functional-roles}

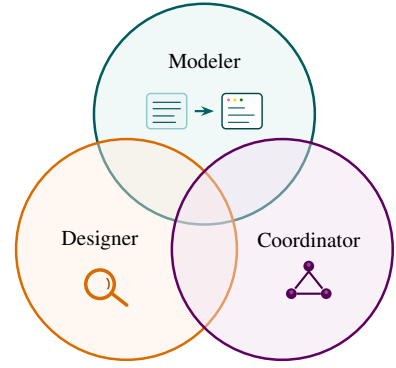
\begin{figure}[!t]
\centering
\begin{tikzpicture}[
  every node/.append style={font=\scriptsize},
  roleText/.style={font=\footnotesize, align=center, text=black}
]
  \draw[draw=teal!72!black, fill=teal!12, fill opacity=0.42,
        line width=0.95pt] (0,0.84) circle [radius=1.46];
  \draw[draw=orange!85!black, fill=orange!13, fill opacity=0.42,
        line width=0.95pt] (-1.03,-0.95) circle [radius=1.46];
  \draw[draw=violet!75!black, fill=violet!12, fill opacity=0.42,
        line width=0.95pt] (1.03,-0.95) circle [radius=1.46];

  \node[roleText] at (0,1.55) {Modeler};
  \begin{scope}[shift={(0,0.88)}, scale=1.06]
    \path[fill=teal!7, draw=teal!44, line width=0.54pt, rounded corners=1.5pt]
      (-0.72,-0.22) rectangle (-0.22, 0.22);
    \draw[teal!64!black, line width=0.44pt, line cap=round]
      (-0.64, 0.11) -- (-0.30, 0.11)
      (-0.64, 0.03) -- (-0.30, 0.03)
      (-0.64,-0.06) -- (-0.38,-0.06)
      (-0.64,-0.14) -- (-0.36,-0.14);
    \draw[teal!78!black, line width=0.82pt, -{Stealth[scale=0.72]}]
      (-0.12, 0) -- (0.12, 0);
    \path[fill=white, draw=teal!55!black, line width=0.54pt, rounded corners=1.5pt]
      (0.22,-0.22) rectangle (0.72, 0.22);
    \fill[red!55]          (0.30, 0.14) circle [radius=0.019];
    \fill[yellow!80!orange](0.38, 0.14) circle [radius=0.019];
    \fill[green!55!black]  (0.46, 0.14) circle [radius=0.019];
    \draw[teal!68!black, line width=0.42pt, line cap=round]
      (0.30, 0.04) -- (0.50, 0.04)
      (0.30,-0.05) -- (0.62,-0.05)
      (0.30,-0.14) -- (0.55,-0.14);
  \end{scope}

  \node[roleText] at (-1.38,-0.82) {Designer};
  \begin{scope}[shift={(-1.38,-1.41)}, scale=1.02]
    \coordinate (mgO) at (0,0);
    \draw[orange!84!black, line width=1.32pt] (mgO) circle[radius=0.185];
    \draw[orange!62!black, line width=0.56pt, line cap=round]
      (40:0.118) arc (40:-30:0.118);
    \draw[orange!84!black, line width=1.32pt, line cap=round]
      ($(mgO)+(-38:0.185)$) -- ++(-38:0.255);
  \end{scope}

  \node[roleText] at (1.38,-0.82) {Coordinator};
  \begin{scope}[shift={(1.38,-1.42)}, scale=0.96]
    \fill[violet!68!black] ( 0.00,  0.26) circle [radius=0.074];
    \fill[violet!68!black] (-0.24, -0.12) circle [radius=0.074];
    \fill[violet!68!black] ( 0.24, -0.12) circle [radius=0.074];
    \draw[violet!56!black, line width=0.80pt]
      ( 0.00,  0.186) -- (-0.188, -0.050)
      ( 0.00,  0.186) -- ( 0.188, -0.050)
      (-0.166, -0.12) -- ( 0.166, -0.12);
    \fill[white, opacity=0.30] ( 0.00,  0.28) circle [radius=0.028];
    \fill[white, opacity=0.30] (-0.24, -0.10) circle [radius=0.028];
    \fill[white, opacity=0.30] ( 0.24, -0.10) circle [radius=0.028];
  \end{scope}

\end{tikzpicture}
\caption{Role classification of LLMs in VRP research.}
\label{fig:llm-capabilities}
\end{figure}

According to how LLMs enter the VRP solving process, existing studies can be grouped into three roles. As modelers, LLMs translate unstructured requirements into variables, constraints, and objective functions. As designers, they generate heuristic code, neighborhood operators, or candidate routes for search frameworks or direct inspection. As coordinators, they manage multi-stage workflows by selecting tools, switching strategies, passing intermediate results, and injecting information into downstream modules. Fig.~\ref{fig:llm-capabilities} illustrates this role classification.

\section{LLMs as Modelers}
\label{sec:modelers}

Transforming delivery requirements into solvable mathematical models is one of the most expertise-intensive stages of the VRP pipeline. Constraints such as time windows, capacities, heterogeneous fleets, driver regulations, service priorities, and depot rules are often scattered across business documents or operational descriptions. Manual translation into variables, objectives, and constraints is laborious, and errors may occur in constraint directions, index ranges, variable semantics, or the correspondence between business rules and mathematical expressions.

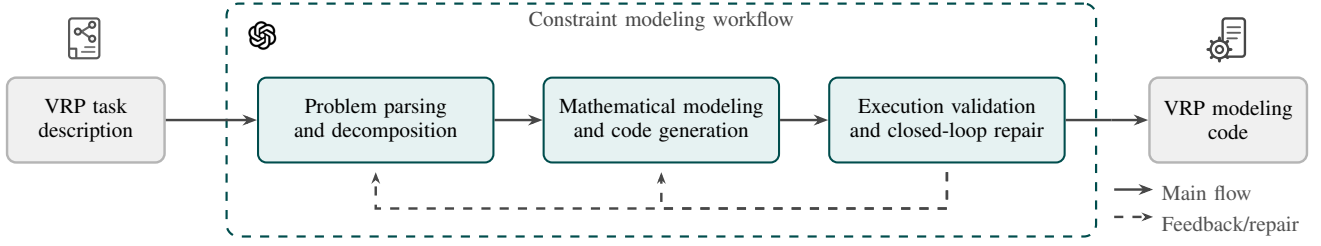
\begin{figure*}[!t]
\centering
\resizebox{0.95\textwidth}{!}{%
\begin{tikzpicture}[
  >=Stealth,
  every node/.append style={font=\scriptsize},
  solvearr/.style={-{Stealth[scale=0.72]}, thick, black!72},
  flowline/.style={thick, black!72},
  feed/.style={-{Stealth[scale=0.58]}, thick, dashed, black!72},
  flowbox/.style={rectangle, rounded corners=2.8pt, minimum width=2.72cm,
                  minimum height=0.98cm, text width=2.38cm,
                  inner xsep=0.12cm, inner ysep=0.10cm,
                  align=center, line width=0.68pt},
  inputbox/.style={rectangle, rounded corners=3pt,
                   minimum width=1.82cm, minimum height=0.98cm,
                   text width=1.58cm,
                   inner xsep=0.10cm, inner ysep=0.10cm,
                   align=center, line width=0.8pt,
                   draw=gray!58, fill=gray!12}
]

  \node[inputbox] (inp) at (0,0) {VRP task\\description};
  \node[inputbox] (out) at (13.20,0) {VRP modeling\\code};
  \taskdescicon{($(inp.north)+(0,0.44)$)}{0.46}
  \modelcodeicon{($(out.north)+(0,0.44)$)}{0.46}

  \node[flowbox, draw=teal!60!black, fill=teal!10] (p1) at (3.35,0)
    {Problem parsing\\and decomposition};

  \node[flowbox, draw=teal!60!black, fill=teal!10] (p2) at (6.65,0)
    {Mathematical modeling\\and code generation};

  \node[flowbox, draw=teal!60!black, fill=teal!10] (p3) at (9.95,0)
    {Execution validation\\and closed-loop repair};

  \node[draw=teal!60!black, dashed, rounded corners=4pt, line width=0.72pt,
        inner xsep=10pt, inner ysep=24pt,
        fit=(p1)(p2)(p3)] (pbox) {};
  \node[anchor=north west, inner sep=0] at ($(pbox.north west)+(0.14,-0.12)$)
    {\includegraphics[height=0.52cm]{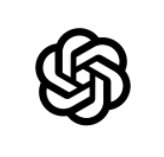}};
  \node[anchor=north, font=\scriptsize, text=black!72,
        fill=white, inner xsep=2pt, inner ysep=0.5pt]
        at ($(pbox.north)+(0,-0.08)$) {Constraint modeling workflow};

  \draw[solvearr] (inp.east) -- (p1.west);
  \draw[solvearr] (p1.east) -- (p2.west);
  \draw[solvearr] (p2.east) -- (p3.west);
  \coordinate (merge) at (11.85,0);
  \draw[flowline] (p3.east) -- (merge);
  \draw[solvearr] (merge) -- (out.west);

  \draw[feed] ($(p3.south)+(0,-0.10)$) -- ++(0,-0.42)
              -| ($(p1.south)+(0,-0.10)$);
  \draw[feed] ($(p3.south)+(0,-0.10)$) -- ++(0,-0.42)
              -| ($(p2.south)+(0,-0.10)$);

  \node[anchor=north east, font=\scriptsize, text=black!68, align=left,
        fill=white, inner sep=1.4pt]
        at ($(out.south east)+(-0.05,-0.18)$) {%
    \raisebox{2pt}{\tikz[baseline=-0.5ex]%
      \draw[solvearr] (0,0)--(0.45,0);}~Main flow\\[2pt]
    \raisebox{3pt}{\tikz[baseline=-0.5ex]%
      \draw[feed] (0,0)--(0.45,0);}~Feedback/repair};
\end{tikzpicture}%
}
\caption{Illustration of the modeler workflow.}
\label{fig:modeler-overview}
\end{figure*}

As modelers, LLMs do not replace the optimizer. Their role is to convert problem descriptions into formal models or modeling code by identifying objects, variables, objectives, and constraint relations in natural language and expressing them in a solver-readable form. Mathematical programming solvers or search algorithms still perform the actual optimization. Fig.~\ref{fig:modeler-overview} summarizes this workflow, which typically includes problem parsing, mathematical modeling, code generation, execution validation, and closed-loop repair.

In general optimization modeling, existing work mainly studies natural-language parsing, mathematical formulation, solver-code generation, and closed-loop repair. An early systematic attempt is OptiMUS 0.2 by AhmadiTeshnizi \textit{et al.} \cite{ahmaditeshnizi2023optimus}, which connects modeling, code generation, execution, verification, and repair to transform natural-language optimization problems into mixed-integer programming (MIP) code. OptiMUS 0.3 \cite{ahmaditeshnizi2024optimus} further introduced problem decomposition and staged validation to improve scalability for large constraint-modeling tasks.

End-to-end modeling workflows form another line of research. Zhang \textit{et al.} \cite{zhang2024optllm} built a pipeline for generating executable optimization models from textual descriptions. Wang \textit{et al.} \cite{wang2025ormind} proposed ORMind, which combines problem understanding, formal modeling, code implementation, and result checking in one reasoning framework. Data and evaluation resources are also important. Huang \textit{et al.} \cite{huang2025orlm} constructed OR-Instruct, proposed ORLM, and evaluated open-source models on realistic operations research scenarios such as IndustryOR. Ma \textit{et al.} \cite{ma2024llamoco} proposed LLaMoCo, which adapts general-purpose LLMs to optimization-code generation through instruction tuning on code-form problem descriptions and expert optimization programs.

Multi-agent and self-improvement mechanisms further extend automatic modeling. Thind \textit{et al.} \cite{thind2025optimai} developed OptimAI with roles such as formulator, planner, coder, and code critic, allowing mathematical formalization, solver-strategy planning, code execution, and correction to proceed collaboratively. Kong \textit{et al.} \cite{kong2025alphaopt} proposed AlphaOPT, where a self-improving experience library stores reusable modeling experience from examples and solver feedback. For problems whose variables, constraints, and objectives are difficult to generate in one pass, Liu \textit{et al.} \cite{liu2025optitree} proposed OptiTree, using a modeling tree and tree search to invoke decomposition at different levels and generate structured prompts. Ding \textit{et al.} \cite{ding2026orr1} combined few-shot SFT with test-time reinforcement learning in OR-R1 to improve data efficiency and output consistency.

For VRP and logistics routing, modeling studies must capture structural relations among vehicles, customers, routes, visit order, and spatiotemporal constraints. The generated model must also connect reliably to solvers and pass feasibility checks, which makes interface consistency and constraint verifiability central. Ju \textit{et al.} \cite{ju2024ttg} proposed TTG for travel-trajectory generation, using SFT to improve output feasibility and emphasizing consistency with exact-solver interfaces. For VRP variants with fine-grained and compound constraints, Jiang \textit{et al.} \cite{jiang2025droc} proposed DRoC, which decomposes composite constraints, retrieves templates from an external knowledge base, and recombines subproblems into modeling code, reducing reliance on memorized model parameters. Yang \textit{et al.} \cite{yang2026orthought} constructed the LogiOR benchmark and proposed ORThought, a structured dual-agent framework that separates abstract mathematical modeling from code implementation and provides a more concrete evaluation setting for logistics modeling.

The modeler pathway remains at an early stage. The most reliable systems usually combine LLM generation with solvers, rule checkers, or template libraries so that formal outputs can be validated and repaired. Important details such as variable references, constraint directions, and index ranges still require external checking. Multi-round repair reduces one-shot errors, but its effectiveness depends on feedback quality. Current studies mainly cover general optimization modeling or a limited set of routing variants; unified modeling frameworks for common problems such as VRPTW and MDVRP remain open.

\section{LLMs as Designers}
\label{sec:designers}

\begin{figure*}[!t]
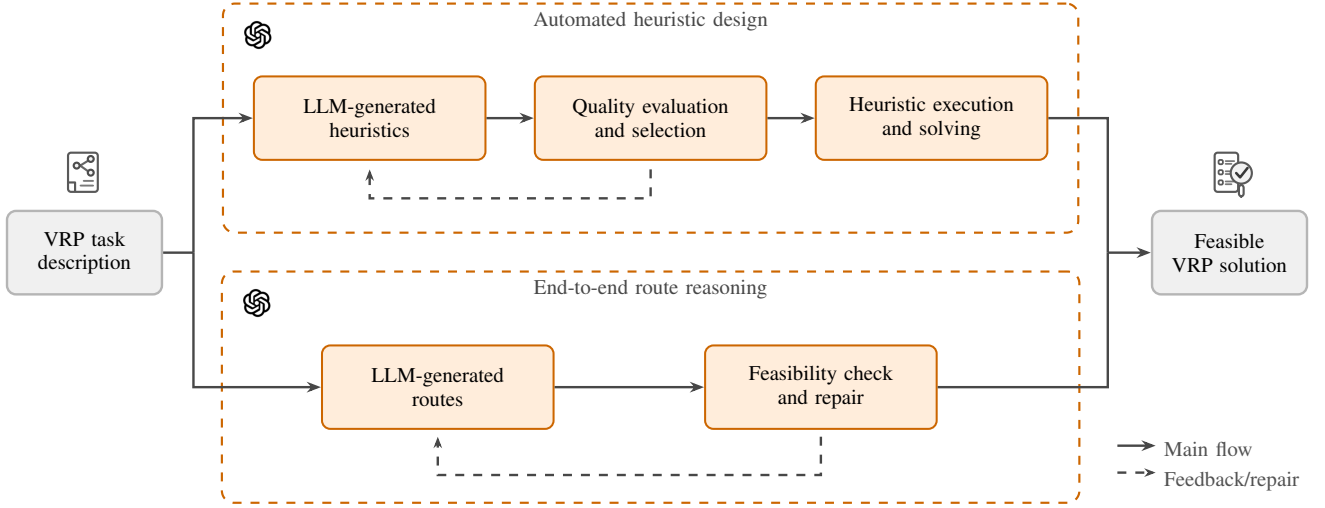

\centering
\resizebox{0.95\textwidth}{!}{%
\begin{tikzpicture}[
  >=Stealth,
  every node/.append style={font={\scriptsize\linespread{1}\selectfont}},
  solvearr/.style={-{Stealth[scale=0.72]}, thick, black!72},
  flowline/.style={thick, black!72},
  feed/.style={-{Stealth[scale=0.58]}, thick, dashed, black!72},
  boxtext/.style={font=\scriptsize,
                  align=center, inner xsep=0.10cm, inner ysep=0.10cm},
  flowbox/.style={rectangle, rounded corners=2.8pt, minimum width=2.72cm,
                  minimum height=0.98cm, text width=2.38cm,
                  boxtext, line width=0.68pt},
  inputbox/.style={rectangle, rounded corners=3pt,
                   minimum width=1.82cm, minimum height=0.98cm,
                   text width=1.58cm,
                   boxtext, line width=0.8pt,
                   draw=gray!58, fill=gray!12}
]
  \node[inputbox] (inp) at (0,0) {VRP task\\description};
  \node[inputbox] (out) at (13.45,0) {Feasible\\VRP solution};
  \taskdescicon{($(inp.north)+(0,0.44)$)}{0.46}
  \solutioncheckicon{($(out.north)+(0,0.44)$)}{0.46}

  \node[flowbox, draw=orange!80!black, fill=orange!14] (h1) at (3.35,1.58)
    {LLM-generated\\heuristics};
  \node[flowbox, draw=orange!80!black, fill=orange!14] (h2) at (6.65,1.58)
    {Quality evaluation\\and selection};
  \node[flowbox, draw=orange!80!black, fill=orange!14] (h3) at (9.95,1.58)
    {Heuristic execution\\and solving};

  \node[draw=orange!80!black, dashed, rounded corners=4pt, line width=0.72pt,
        inner xsep=10pt, inner ysep=24pt,
        fit=(h1)(h2)(h3)] (hbox) {};
  \node[anchor=north west, inner sep=0] at ($(hbox.north west)+(0.14,-0.12)$)
    {\includegraphics[height=0.52cm]{figs/LLM_image.png}};
  \node[anchor=north, font=\scriptsize, text=black!72,
        fill=white, inner xsep=2pt, inner ysep=0.5pt]
        at ($(hbox.north)+(0,-0.08)$) {Automated heuristic design};

  \node[flowbox, draw=orange!80!black, fill=orange!14] (e1) at (4.15,-1.58)
    {LLM-generated\\routes};
  \node[flowbox, draw=orange!80!black, fill=orange!14] (e2) at (8.65,-1.58)
    {Feasibility check\\and repair};
  \coordinate (eboxnw) at ($(hbox.north west)+(0,-3.16)$);
  \coordinate (eboxse) at ($(hbox.south east)+(0,-3.16)$);

  \node[draw=orange!80!black, dashed, rounded corners=4pt, line width=0.72pt,
        inner sep=0pt,
        fit=(eboxnw)(eboxse)] (ebox) {};
  \node[anchor=north west, inner sep=0] at ($(ebox.north west)+(0.14,-0.12)$)
    {\includegraphics[height=0.52cm]{figs/LLM_image.png}};
  \node[anchor=north, font=\scriptsize, text=black!72,
        fill=white, inner xsep=2pt, inner ysep=0.5pt]
        at ($(ebox.north)+(0,-0.08)$) {End-to-end route reasoning};

  \coordinate (split) at (1.28,0);
  \coordinate (merge) at (12.02,0);
  \draw[flowline] (inp.east) -- (split);
  \draw[solvearr] (split) |- (h1.west);
  \draw[solvearr] (split) |- (e1.west);

  \draw[solvearr] (h1.east) -- (h2.west);
  \draw[solvearr] (h2.east) -- (h3.west);
  \draw[flowline] (h3.east) -- ++(0.42,0) -- ($(merge)+(0,1.58)$) -- (merge);

  \draw[solvearr] (e1.east) -- (e2.west);
  \draw[flowline] (e2.east) -- ++(0.42,0) -- ($(merge)+(0,-1.58)$) -- (merge);
  \draw[solvearr] (merge) -- (out.west);

  \draw[feed] ($(h2.south)+(0,-0.10)$) -- ++(0,-0.34)
              -| ($(h1.south)+(0,-0.10)$);
  \draw[feed] ($(e2.south)+(0,-0.08)$) -- ++(0,-0.45)
              -| ($(e1.south)+(0,-0.08)$);

  \node[anchor=south east, font=\scriptsize, text=black!68, align=left,
        fill=white, inner sep=1.4pt]
        at ($(out.south east)+(-0.05,-2.35)$) {%
    \raisebox{2pt}{\tikz[baseline=-0.5ex]\draw[solvearr] (0,0)--(0.45,0);}~Main flow\\[2pt]
    \raisebox{3pt}{\tikz[baseline=-0.5ex]\draw[feed] (0,0)--(0.45,0);}~Feedback/repair};
\end{tikzpicture}%
}
\caption{Illustration of the designer workflow.}
\label{fig:designer-workflow}
\end{figure*}

As designers, LLMs generate either solving components or candidate solutions. One line of work asks the model to produce heuristic algorithms, neighborhood operators, or executable programs, which are then evaluated repeatedly by an external evaluator. Another line bypasses explicit solvers and asks the model to generate route sequences directly. Fig.~\ref{fig:designer-workflow} summarizes these two workflows.

\subsection{Automated Heuristic Design}
\label{subsec:heuristic-design}

Automated heuristic design embeds LLMs into evolutionary or iterative search loops. The model generates or rewrites heuristic functions, and an external evaluator scores, filters, and retains useful candidates. Solution quality depends less on a single model output than on repeated execution feedback and selection pressure.

Early studies centered on candidate generation and external evaluation. Liu \textit{et al.} \cite{liu2023algorithm} introduced LLMs into evolutionary algorithms through AEL, where new candidate functions are generated from existing heuristic code. Romera-Paredes \textit{et al.} \cite{romeraparedes2024mathematical} proposed FunSearch, in which LLMs generate executable program fragments and an external scoring function selects programs that improve the target metric. Liu \textit{et al.} \cite{liu2024eoh} further proposed EoH, connecting heuristic generation, execution evaluation, population update, and candidate elimination into a closed loop for automatic algorithm design.

Later studies introduced reflection, tree search, and diversity-preserving mechanisms. Ye \textit{et al.} \cite{ye2024reevo} proposed ReEvo, where the LLM first diagnoses weak candidates and then generates improved versions. Zheng \textit{et al.} \cite{zheng2025mctsahd} introduced Monte Carlo tree search through MCTS-AHD, using a tree to organize and expand heuristic candidates. Dat \textit{et al.} \cite{dat2025hsevo} proposed HSEvo by combining harmony search and genetic algorithms to preserve diversity and reduce premature convergence. Yao \textit{et al.} \cite{yao2025meoh} proposed MEoH, which jointly considers solution quality and heuristic diversity so that elites can improve without collapsing the candidate pool.

In VRP, the emphasis has shifted from generic heuristic generation to use within concrete solver frameworks. Ma \textit{et al.} \cite{ma2024autodh} represented candidate heuristics as executable code and used VRP feedback to select improved versions, showing that LLM-generated heuristics can support routing. Hottung \textit{et al.} \cite{hottung2025vrpagent} proposed VRPAgent, moving from individual heuristic functions toward automatic discovery, evaluation, and retention of heuristic operators. Malik \textit{et al.} \cite{malik2026pyvrpplus} proposed PyVRP+, which builds on the HGS framework in PyVRP and lets LLMs improve local components through reasoning, execution, and reflection. Xie \textit{et al.} \cite{xie2025ailsahd} similarly generated replaceable or complementary local heuristic components around an existing CVRP solver, reflecting a move from standalone heuristic generation to plug-in-style solver enhancement.

Some work constrains the generation space more explicitly. Zhao \textit{et al.} \cite{zhao2025glns} proposed G-LNS, requiring LLMs to generate destroy and repair operators for LNS. This confines outputs to key search operations, reduces invalid programs, and makes the generated operators easier to use with existing procedures. Shi \textit{et al.} \cite{shi2026llmvd} extended automated heuristic design to vehicle-drone routing, using structured prompts and self-debugging to handle vehicle-drone coordination constraints.

Recent work also studies framework-level closed-loop design and learning-enhanced generation. Liu \textit{et al.} \cite{liu2024llm4ad} released LLM4AD, a platform that standardizes prompts, code execution, result evaluation, and log management for reproducible algorithm-design experiments. Yang \textit{et al.} \cite{yang2025heurageni} proposed HeurAgenix, where LLM agents coordinate heuristic generation, evaluation, and selection for complex combinatorial optimization tasks. Shi \textit{et al.} \cite{shi2026moh} used meta-optimization to generate more generalizable heuristics, extending the search target from a single heuristic function to a heuristic optimizer. Huang \textit{et al.} \cite{huang2025calm} proposed CALM, using co-evolutionary fine-tuning between algorithms and LLMs to bias generation toward more effective expressions. Zhu \textit{et al.} \cite{zhu2025rfthgs} fine-tuned LLMs with reinforcement learning to improve key HGS components for CVRP. Cao \textit{et al.} \cite{cao2025llmqlearningcvrptw} proposed an LLM-enhanced Q-learning method for CVRPTW, where LLMs assist exploration, self-checking, and prioritized replay. These studies suggest that LLMs are evolving from program generators into strategy-assistance modules within learning processes.

Automated heuristic design shows a clear pattern: LLMs propose programs, while external fitness evaluation provides quality control. Simple random mutation is gradually being replaced by reflection, diversity preservation, multi-objective selection, and execution-feedback-driven refinement. The limitations are also clear. Frequent LLM calls and external evaluation introduce overhead, and the stability, transferability, and reproducibility of generated heuristics remain sensitive to problem interfaces, scoring functions, and execution environments.

\subsection{End-to-End Route Reasoning}
\label{subsec:route-reasoning}

End-to-end route reasoning feeds a natural-language VRP instance directly into an LLM and asks it to output customer visiting sequences. The solving process is largely internal to the model; external programs or solvers appear only in extended settings. Huang \textit{et al.} \cite{huang2024words} systematically evaluated the direct solving ability of models such as GPT-4 on standard VRP benchmarks. The results show that LLMs can produce reasonable routes for small instances. With chain-of-thought prompting, the model can track capacity constraints and route validity to some extent, showing limited zero-shot adaptability without additional training.

Li \textit{et al.} \cite{li2025ars} proposed ARS, further exploring how LLMs can assemble automated routing-solving pipelines by combining natural-language understanding with automatic organization of routing procedures. For specific variants, Zafar \textit{et al.} \cite{zafar2025evrp} studied advanced prompting and self-checking mechanisms for EVRP. The results suggest that clearer prompt design and self-verification can improve the ability of the model to handle additional constraints without extra training.

End-to-end route reasoning has low deployment cost and requires neither training nor complex external interfaces. It may also adapt quickly to new variants described in natural language. However, as the number of customers grows, context limits and combinatorial reasoning weaknesses become severe, and the quality gap from specialized exact algorithms and metaheuristics widens rapidly. Most such methods also lack formal feasibility guarantees, which makes direct deployment risky when hard constraints such as time windows and vehicle capacities must be respected.

\section{LLMs as Coordinators}
\label{sec:coordinators}

\begin{figure*}[!t]
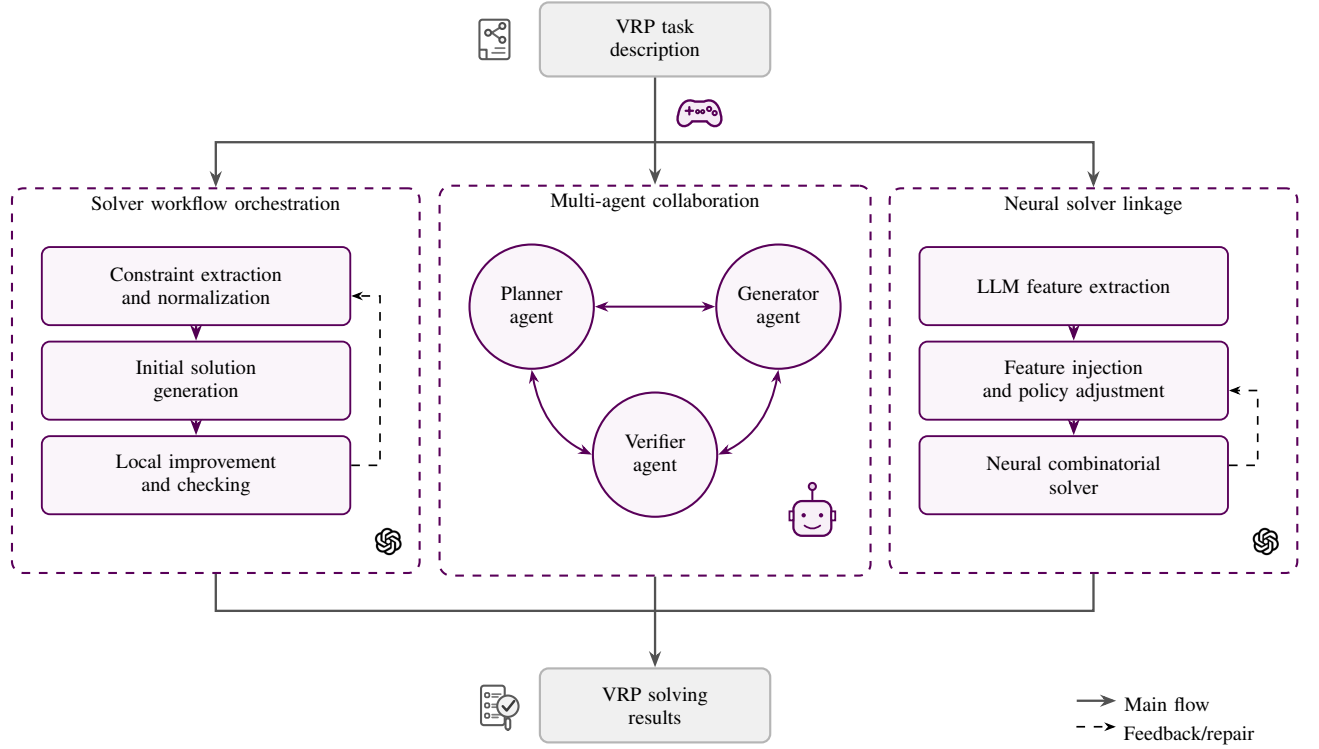

\centering
\resizebox{0.96\textwidth}{!}{%
\begin{tikzpicture}[
  x=1cm, y=1cm,
  >=Stealth,
  every node/.append style={font=\scriptsize},
  boxtext/.style={font=\scriptsize,
                  align=center, inner xsep=0.10cm, inner ysep=0.10cm},
  iobox/.style={rectangle, rounded corners=3pt, minimum width=2.8cm, minimum height=0.90cm,
                text width=2.45cm,
                draw=gray!58, fill=gray!12, boxtext, line width=0.75pt},
  stepA/.style={rectangle, rounded corners=2.8pt, minimum width=3.75cm, minimum height=0.95cm,
                text width=3.45cm, boxtext,
                draw=violet!70!black, fill=violet!4,
                line width=0.65pt},
  agboxB/.style={circle, draw=violet!70!black, fill=violet!4, line width=0.72pt,
                 minimum size=1.45cm, text width=1.14cm, boxtext},
  stepC/.style={rectangle, rounded corners=2.8pt, minimum width=3.75cm, minimum height=0.95cm,
                text width=3.45cm, boxtext,
                draw=violet!70!black, fill=violet!4,
                line width=0.65pt},
  arr/.style={-{Stealth[scale=0.85]}, line width=0.80pt, black!68},
  flowline/.style={line width=0.80pt, black!68},
  arrdash/.style={-{Stealth[scale=0.72]}, line width=0.58pt, dashed, black},
  arrbothB/.style={<->, line width=0.70pt, violet!70!black,
                   arrows={Stealth[scale=0.76]-Stealth[scale=0.76]}}
]

\path[use as bounding box] (-0.5, 0.75) rectangle (15.5, -9.10);

\node[stepA] (t1) at (2.03, -3.15) {Constraint extraction\\and normalization};
\node[stepA] (t2) at (2.03, -4.30) {Initial solution\\generation};
\node[stepA] (t3) at (2.03, -5.45) {Local improvement\\and checking};
\draw[arr, violet!65!black] (t1) -- (t2);
\draw[arr, violet!65!black] (t2) -- (t3);
\draw[arrdash] ($(t3.east)+(0,0.12)$) -- ++(0.35,0) |- ($(t1.east)+(0,-0.12)$);

\coordinate (box1-east-aid) at ($(t1.east)+(0.48,0)$);
\node[draw=violet!70!black, dashed, rounded corners=4pt, line width=0.72pt,
      inner xsep=10pt, inner ysep=20pt,
      fit=(t1)(t2)(t3)(box1-east-aid)] (box1) {};
\node[anchor=north, font=\scriptsize, text=black,
      fill=white, inner xsep=2pt, inner ysep=0.5pt]
    at ($(box1.north)+(0,-0.08)$) {Solver workflow orchestration};
\node[anchor=south east, inner sep=0] at ($(box1.south east)+(-0.12,0.12)$)
  {\includegraphics[height=0.52cm]{figs/LLM_image.png}};

\node[agboxB] (a1) at (6.11, -3.40) {Planner\\agent};
\node[agboxB] (a2) at (9.11, -3.40) {Generator\\agent};
\node[agboxB] (a3) at (7.61, -5.20) {Verifier\\agent};
\draw[arrbothB] (a1.east) -- (a2.west);
\draw[arrbothB] (a1.south) to[out=280, in=155] (a3.west);
\draw[arrbothB] (a2.south) to[out=260, in=25] (a3.east);

\node[draw=violet!70!black, dashed, rounded corners=4pt, line width=0.72pt,
      inner xsep=10pt, inner ysep=20pt,
      fit=(a1)(a2)(a3)] (box2) {};
\node[anchor=north, font=\scriptsize, text=black,
      fill=white, inner xsep=2pt, inner ysep=0.5pt]
    at ($(box2.north)+(0,-0.08)$) {Multi-agent collaboration};
\begin{scope}[shift={($(box2.south east)+(-0.72,0.66)$)}, scale=0.48]
  \draw[violet!68!black, line width=0.70pt, line cap=round]
    (0,0.58) -- (0,0.84);
  \draw[violet!68!black, fill=violet!18, line width=0.70pt]
    (0,0.91) circle (0.09);
  \draw[violet!68!black, fill=violet!6, line width=0.78pt, rounded corners=3pt]
    (-0.48,-0.34) rectangle (0.48,0.54);
  \draw[violet!68!black, fill=white, line width=0.62pt]
    (-0.58,0.05) rectangle (-0.48,0.28)
    (0.48,0.05) rectangle (0.58,0.28);
  \fill[violet!70!black] (-0.22,0.18) circle (0.055);
  \fill[violet!70!black] (0.22,0.18) circle (0.055);
  \draw[violet!70!black, line width=0.62pt, line cap=round]
    (-0.20,-0.08) .. controls (-0.06,-0.18) and (0.06,-0.18) .. (0.20,-0.08);
\end{scope}

\node[stepC] (n1) at (12.7, -3.15) {LLM feature extraction};
\node[stepC] (n2) at (12.7, -4.30) {Feature injection\\and policy adjustment};
\node[stepC] (n3) at (12.7, -5.45) {Neural combinatorial\\solver};
\draw[arr, violet!65!black] (n1) -- (n2);
\draw[arr, violet!65!black] (n2) -- (n3);
\draw[arrdash] ($(n3.east)+(0,0.12)$) -- ++(0.35,0) |- ($(n2.east)+(0,-0.12)$);

\coordinate (box3-east-aid) at ($(n1.east)+(0.48,0)$);
\node[draw=violet!70!black, dashed, rounded corners=4pt, line width=0.72pt,
      inner xsep=10pt, inner ysep=20pt,
      fit=(n1)(n2)(n3)(box3-east-aid)] (box3) {};
\node[anchor=north, font=\scriptsize, text=black,
      fill=white, inner xsep=2pt, inner ysep=0.5pt]
    at ($(box3.north)+(0,-0.08)$) {Neural solver linkage};
\node[anchor=south east, inner sep=0] at ($(box3.south east)+(-0.12,0.12)$)
  {\includegraphics[height=0.52cm]{figs/LLM_image.png}};

\coordinate (tbar) at (box2.north |- 0,-1.40);
\node[iobox] (ctrl) at ($(tbar)+(0,1.25)$) {VRP task\\description};
\taskdescicon{($(ctrl.west)+(-0.52,0)$)}{0.48}
\draw[flowline] (ctrl.south) -- (tbar);
\controllericon{($(ctrl.south)!0.50!(tbar)+(0.54,0)$)}{0.56}

\draw[arr] (tbar) -| (box1.north);
\draw[arr] (tbar) -| (box2.north);
\draw[arr] (tbar) -| (box3.north);

\coordinate (bbar) at (box2.south |- 0,-7.10);
\draw[flowline] (box1.south) |- (bbar);
\draw[flowline] (box2.south) |- (bbar);
\draw[flowline] (box3.south) |- (bbar);

\node[iobox] (out) at ($(bbar)+(0,-1.15)$) {VRP solving\\results};
\solutioncheckicon{($(out.west)+(-0.52,0)$)}{0.48}
\draw[arr] (bbar) -- (out.north);

\node[anchor=south east, font=\scriptsize, text=black, align=left,
      fill=white, inner sep=1.4pt]
  at ($(15.5,-9.10)+(-0.55,0.32)$) {%
  \raisebox{2pt}{\tikz[baseline=-0.5ex,>=Stealth]%
    \draw[-{Stealth[scale=0.85]}, line width=0.80pt, black!68] (0,0)--(0.48,0);}~Main flow\\[2pt]
  \raisebox{3pt}{\tikz[baseline=-0.5ex,>=Stealth]%
    \draw[-{Stealth[scale=0.72]}, line width=0.58pt, dashed, black] (0,0)--(0.48,0);}~Feedback/repair};
\end{tikzpicture}%
}
\caption{Illustration of the coordinator workflow.}
\label{fig:coordinator-overview}
\end{figure*}

The coordinator role concerns multi-stage solving rather than isolated model generation or heuristic design. In this role, LLMs decompose tasks, invoke tools, switch modules, and pass intermediate results across a routing pipeline. As illustrated in Fig.~\ref{fig:coordinator-overview}, three pathways are most common: solver workflow orchestration, multi-agent collaboration, and connections with neural solvers.

\subsection{Solver Workflow Orchestration}
\label{subsec:workflow-orchestration}

Solver workflow orchestration lets LLMs take upper-level control in multi-stage solving processes. This line is closely related to LLM-agent mechanisms that interleave reasoning, action execution, and feedback reading. ReAct \cite{yao2022react} showed that LLMs can alternate between reasoning traces and external actions while preserving intermediate states. Toolformer \cite{schick2023toolformer} showed that models can learn when to call tools and how to use returned results during generation. AutoGPT \cite{autogpt} and LangChain \cite{langchain} provide widely used engineering frameworks for coordinating memory, planning, executors, and external APIs.

In VRP, these mechanisms appear when a coordinator organizes the pipeline according to the current instance state. A complex delivery task may require constraint extraction, data normalization, initial-solution construction, local improvement, feasibility checking, and result rewriting. If different modules handle these stages, the LLM acts as a process controller: it selects the next tool or submodule and converts previous outputs into representations that the next stage can consume. The quality of this process depends strongly on interface consistency and execution order.

In broader combinatorial optimization, Chen \textit{et al.} \cite{chen2026dragon} proposed DRAGON, decomposing large-scale problems into global solution analysis, local subproblem reconstruction, and feedback-memory update. LLM agents identify promising improvement regions and organize local optimization. For VRP specifically, however, systematically evaluated workflow-orchestration mechanisms remain limited. Most studies still adapt general agent frameworks to routing tasks. This line of work suggests that LLMs can assume workflow-level control, but stable performance gains on concrete VRP benchmarks still need stronger evidence.

\subsection{Multi-Agent Collaboration}
\label{subsec:multi-agent-collaboration}

Multi-agent collaboration decomposes complex solving into several LLM agents with explicit roles. In general combinatorial optimization, Wang \textit{et al.} \cite{wang2025maef} proposed MAEF, decomposing scheduling optimization into problem definition, initial population generation, evolutionary optimization, and result evaluation. Yan \textit{et al.} \cite{yan2026hmace} proposed HMACE, organizing heuristic search as collaborative evolution among proposer, generator, evaluator, and reflector agents; experiments on TSP, bin packing, and scheduling suggest that role specialization can improve search efficiency and solution quality. Qu \textit{et al.} \cite{qu2026coral} proposed CORAL for open-ended discovery, combining shared persistent memory, asynchronous agent execution, and candidate evaluation to support long-horizon search and knowledge accumulation.

In existing VRP research, multi-agent collaboration is one of the most direct manifestations of the coordinator role. Zhang \textit{et al.} \cite{zhang2025afl} proposed AFL, a multi-agent framework for end-to-end automated solving of complex VRPs. AFL decomposes the overall task into manageable subtasks and assigns them to agents with different functions. The main idea is to reduce error accumulation in long-chain reasoning through role separation and to improve consistency and feasibility through cross-checking among agents.

Multi-agent collaboration reduces the fragility of a single long reasoning chain, but it also introduces communication overhead and new paths for error propagation. If an intermediate agent makes a systematic mistake, downstream modules may continue under a wrong premise. Role boundaries, interaction protocols, and stopping criteria are often set empirically. Multi-agent collaboration is therefore better viewed as a system paradigm for making complex workflows more manageable than as a mature standard solution.

\subsection{Connections with Neural Solvers}
\label{subsec:neural-solver-connections}

Another coordinator pathway combines LLMs with downstream neural solvers or learned policy networks. In this setting, LLMs act as high-level strategy controllers. They provide instance-level guidance, inject biases, adjust policies, or generate constraint masks for existing neural solvers, thereby influencing search directions, attention distributions, or feasible action sets.

Existing studies can be grouped into three categories. The first injects semantic or structural information obtained by LLMs into neural solvers, allowing the original solving backend to receive instance-level guidance. Tran \textit{et al.} \cite{tran2025llmnco} proposed NCO-LLMs, which retain an NCO model as the execution backend. The LLM extracts structural features from an instance, converts them into guiding biases in attention mechanisms, and improves transfer to large-scale instances through fine-tuning on small but diverse data. Malik \textit{et al.} \cite{malik2026llmaide} proposed LLMAide with a similar motivation. It fuses semantic features extracted by LLMs with spatial representations in neural routing solvers through a multi-scale fusion module and validates the effectiveness of this fusion on multiple VRP variants.

The second category targets unified handling of cross-variant constraints. Zhou \textit{et al.} \cite{zhou2025urs} proposed URS, a unified neural routing solver that improves zero-shot generalization through unified data representation and a parameter generator. It converts natural-language constraint descriptions into stepwise mask functions, allowing complex constraints to enter neural construction. Chi \textit{et al.} \cite{chi2026unslhe} proposed a unified framework that uses LLM evolution to generate variant-specific heuristics and adjust neural-solver attention distributions, reducing the need to retrain for every variant.

The third category starts from representation alignment and task transfer. Zeng \textit{et al.} \cite{zeng2026universalvrp} encoded problem descriptions as task embeddings and fed them together with node embeddings into an encoder-decoder structure. Pretraining, LoRA fine-tuning, and mixture-of-experts mechanisms were used to improve transfer across routing problems such as TSP and CVRP. Feng \textit{et al.} \cite{feng2026alignopt} proposed AlignOPT for general combinatorial optimization, aligning semantic representations from LLMs with structural representations from graph neural solvers to improve the generalization of neural heuristics across instances.

The coordinator pathway is developing quickly, although empirical evidence remains scattered. Existing studies show that LLMs can move beyond modeling and design to control multi-module systems through task decomposition, tool invocation, and solver coordination. Connections with neural solvers can be viewed as adding high-level cognitive control on top of mature learned backends. This category should be limited to cases where LLMs explicitly participate in input construction, feature guidance, policy correction, constraint-mask generation, or heuristic evolution for neural solvers, rather than all neural routing models. Public evaluation resources, unified protocols, and reproducible experiments for coordinator-style methods remain insufficient, and existing studies still differ substantially in system boundaries and evaluation goals.

\section{Datasets and Evaluation Resources}
\label{sec:datasets}

Datasets and evaluation resources determine whether methods can be compared under consistent conditions. For LLM-assisted VRP solving, an instance set alone is not sufficient. A complete evaluation should specify how the problem is verbalized, which tools the model may call, how outputs are parsed, and how solution quality, feasibility, and computational cost are measured. Existing resources can be grouped by data source and evaluation goal into standardized VRP benchmarks, real or near-real operational datasets, and LLM-oriented evaluation frameworks. These resources support classical algorithm comparison, real-scenario generalization analysis, and model-capability assessment, respectively.

Table~\ref{tab:dataset_resource_summary} summarizes the standardized and operational datasets, and Table~\ref{tab:llm_eval_frameworks} lists the LLM-related evaluation frameworks. The category boundaries are not strict: some LLM-oriented tasks still use instances from CVRPLIB, Solomon, or operational data. The additional value of these resources lies in textual descriptions, interaction interfaces, executable environments, feedback mechanisms, and evaluation protocols.

\begin{table*}[!t]
\centering
\caption{Summary of VRP-related datasets.}
\label{tab:dataset_resource_summary}
\footnotesize
\setlength{\tabcolsep}{3pt}
\renewcommand{\arraystretch}{1.14}
\begin{tabular*}{\textwidth}{@{\extracolsep{\fill}}>{\raggedright\arraybackslash}p{0.18\textwidth}>{\raggedright\arraybackslash}p{0.21\textwidth}>{\raggedright\arraybackslash}p{0.55\textwidth}@{}}
\toprule
Dataset & Covered tasks & Scale and data characteristics \\
\midrule
\rowcolor{tableTitleColor}
\multicolumn{3}{c}{Standardized benchmark datasets} \\
CVRPLIB \cite{uchoa2017new} & CVRP with capacity constraints & A widely used collection of 376 CVRP instances with 12--30000 customers,
 covering different scales, customer distributions, demand patterns, and capacity settings. \\
Solomon \cite{solomon1987algorithms} & VRPTW with time windows & 56 instances with 100 customers, grouped by random, clustered,
 and mixed customer distributions and different route-time structures. \\
Archetti \textit{et al.} \cite{archetti2006tabu} & SDVRP with split deliveries & Instances derived from Christofides cases 1--5, 11, 
and 12, with about 50--199 customers and demand that may be split across vehicles. \\
Schneider \textit{et al.} \cite{schneider2014electric} & E-VRPTW with charging and battery constraints & Solomon-style electric VRP instances, 
including 36 small cases with 5--15 customers and 2--8 charging stations, 56 large cases with 100 customers and 21 charging stations. \\
Cordeau \textit{et al.} \cite{cordeau1997tabu} & MDVRP with multiple depots & 33 multi-depot instances with about 48--360 customers and 2--9 depots,
 involving depot selection, customer assignment, vehicle capacity, and route-duration constraints. \\
Goetschalckx \textit{et al.} \cite{goetschalckx1989vehicle} & VRPB with backhauls & 62 instances with about 20--150 customers, 
distinguishing linehaul and backhaul customers and requiring deliveries before pickups. \\
Li \textit{et al.} \cite{li2007open} & OVRP with open routes & 16 medium-scale instances with 50--199 customers and 8 large-scale instances with 200--480 customers,
 where vehicles do not return to the depot after service. \\
Christofides \textit{et al.} \cite{christofides1979vehicle} & VRP with route-length limits & 14 classical instances with 50--199 customers;
 cases C6--C10 and C13--C14 include maximum route-length constraints. \\
\midrule
\rowcolor{tableTitleColor}
\multicolumn{3}{c}{Real or near-real operational datasets} \\
Amazon Last-Mile \cite{merchan2022amazonlastmile} & Real last-mile delivery analysis & 9184 historical routes from 2018, 
including 6112 training routes and 3072 evaluation routes from five U.S. metropolitan areas with route, stop, and package-level features. \\
Delivering Data \cite{vrani2025deliveringdata} & Real pharmaceutical last-mile delivery & Nine days of third-party pharmaceutical delivery data,
 with about 60--85 delivery points per day and distance and travel-time matrices under different traffic scenarios. \\
Olist-VRP \cite{greenberg2025earli} & Near-real CVRP construction & 4096 CVRP instances with 20--500 customers, constructed from Brazilian e-commerce orders, seller locations, and urban travel-time information. \\
SVRPBench \cite{heakl2025svrpbench} & Near-real stochastic and dynamic VRP & More than 500 instances with 10--1000 customers, 
simulating time-dependent congestion, stochastic delays, accidents, and empirical time windows. \\
EvoReal \cite{zhu2026evoreal} & Near-real routing-instance generation & LLM-guided generation of realistic synthetic TSP/CVRP training instances to reduce the distribution gap between synthetic training data and real benchmarks. \\
\bottomrule
\end{tabular*}
\end{table*}

\subsection{Standardized Benchmark Datasets}

Standardized VRP benchmarks have long been used for algorithm comparison. CVRPLIB and its related instance families \cite{uchoa2017new} are among the most frequently used resources for the capacitated VRP and support evaluation of solution quality, runtime, optimality gap, and scalability. For the VRP with time windows, the Solomon benchmark \cite{solomon1987algorithms} constructs instances with random, clustered, and mixed customer distributions. It is widely used to compare feasible-route construction and route optimization under service-time constraints.

Several benchmark families correspond to more specific VRP variants. For SDVRP, Archetti \textit{et al.} \cite{archetti2006tabu} constructed instances in which the demand of one customer may be divided among multiple vehicles. Schneider \textit{et al.} \cite{schneider2014electric} proposed E-VRPTW instances that combine customer visits, charging stations, battery consumption, and time-window constraints, making them common references for electric routing studies. For multi-depot routing, Cordeau \textit{et al.} \cite{cordeau1997tabu} provided instances involving depot assignment and route construction. VRPB is often evaluated using the instances of Goetschalckx \textit{et al.} \cite{goetschalckx1989vehicle}, where delivery and pickup customers must be ordered according to linehaul-backhaul rules. OVRP studies commonly use the test problems of Li, Golden, and Wasil \cite{li2007open}, while the classical instances of Christofides \textit{et al.} \cite{christofides1979vehicle} provide references for route-length or route-duration constraints.

Standardized benchmarks have clear problem definitions, relatively complete best-known or optimal solutions, and a long record of academic use, making them suitable for controlled comparison across algorithms. Their limitations are also clear. Many instances are based on stylized coordinate distributions, static demands, and fixed objective functions, and they differ from practical delivery settings that involve textual constraints, time-varying traffic, stochastic disruptions, and dynamic orders.

\begin{table*}[!t]
\centering
\caption{Summary of LLM-related evaluation frameworks.}
\label{tab:llm_eval_frameworks}
\footnotesize
\setlength{\tabcolsep}{3pt}
\renewcommand{\arraystretch}{1.14}
\begin{tabular*}{\textwidth}{@{\extracolsep{\fill}}>{\raggedright\arraybackslash}p{0.18\textwidth}>{\raggedright\arraybackslash}p{0.21\textwidth}>{\raggedright\arraybackslash}p{0.55\textwidth}@{}}
\toprule
Framework & Evaluation target & Scale and characteristics \\
\midrule
RoutBench \cite{li2025ars} & Automatic routing-solver evaluation & 1000 VRP tasks constructed from 24 combinations of routing attributes and constraints, each with problem descriptions, instance data, and feasibility checks. \\
HeuriGym \cite{chen2025heurigym} & Heuristic generation with execution feedback & 
Covers 9 combinatorial optimization problems from 4 domains and 218 instances; the logistics tasks include crew pairing and pickup-and-delivery with time windows. \\
RL4CO \cite{berto2025rl4co} & Learning-based CO training and evaluation & 
Provides 27 combinatorial optimization environments and 23 reference methods; routing environments can generate TSP, CVRP, PDP, MTVRP, and related instances. \\
CO-Bench \cite{sun2025cobench} & Combinatorial-optimization algorithm search & 
An evaluation suite of 36 real-world CO problems across 8 categories, including TSP and PVRP, for testing LLM agents that design algorithms. \\
\bottomrule
\end{tabular*}
\end{table*}

\subsection{Real or Near-Real Operational Datasets}

Real and near-real operational datasets are closer to practical delivery processes and usually come from enterprise operation logs, e-commerce orders, or high-fidelity simulation. The Amazon Last-Mile dataset \cite{merchan2022amazonlastmile} provides historical last-mile routes, stops, and package-level features from five U.S. metropolitan areas. It is useful for analyzing executed route patterns, stop-order preferences, and generalization under real operational distributions. Vrani \textit{et al.} \cite{vrani2025deliveringdata} released a pharmaceutical last-mile delivery dataset from a third-party logistics provider, including multi-day delivery tasks, distance matrices, and travel-time matrices under different traffic scenarios.

Recent work has also constructed near-real VRP benchmarks from public order data. Olist-VRP \cite{greenberg2025earli} uses Brazilian e-commerce orders, seller locations, and road travel times to generate CVRP instances closer to urban delivery settings. This design supports evaluation of learning-based methods under nonuniform spatial distributions and realistic distance structures. SVRPBench \cite{heakl2025svrpbench} further introduces time-dependent congestion, stochastic delays, accident disruptions, and empirical time windows, providing a higher-fidelity benchmark for stochastic and dynamic routing. Instead of directly releasing fixed operational instances, EvoReal \cite{zhu2026evoreal} addresses transfer from synthetic training data to real routing distributions by using LLMs to guide the generation of realistic training instances. This approach reduces the generalization gap of neural solvers on TSPLIB- and CVRPLIB-style benchmarks.

Real or near-real data can better reflect customer spatial distributions, traffic-time variation, historical route regularities, and demand fluctuation. They are valuable for testing robustness and generalization of LLM-assisted or learning-based methods. However, optimal or best-known solutions are usually less complete than in CVRPLIB-style benchmarks. Even when reference solutions are available, objective definitions, anonymization procedures, and business rules are often scenario-dependent. Such datasets are best used as complementary validation rather than replacements for standardized benchmarks.

\subsection{LLM-Oriented Evaluation Frameworks}

LLM-oriented evaluation resources focus more on task organization and model-capability assessment than on numerical instances alone. They typically provide task interfaces, environment wrappers, baseline methods, feedback signals, and relatively unified protocols. RoutBench \cite{li2025ars} is especially relevant to LLM-assisted routing. It constructs 1000 VRP variants from 24 combinations of routing attributes and constraints and evaluates automatic routing solvers on their ability to interpret complex constraints and generate feasible solutions.

RL4CO \cite{berto2025rl4co} is closer to a unified experimental framework for learning-based solvers. Although originally designed for NCO and reinforcement learning, it provides standardized environments, training configurations, and baselines, and can randomly generate multiple routing problem types such as TSP, CVRP, PDP, and MTVRP. HeuriGym \cite{chen2025heurigym} focuses on evaluating LLM-generated heuristics. It emphasizes tool use, executable code feedback, and heuristic quality, making it more suitable for method-level evaluation of designer-style LLM systems. CO-Bench \cite{sun2025cobench} further evaluates LLM agents from the perspective of combinatorial-optimization algorithm search, covering multiple real-world optimization tasks and including routing-related problems such as TSP and PVRP.

These frameworks help identify the capability boundaries of LLMs in complex constraint understanding, heuristic code generation, external solver invocation, cross-variant adaptation, and iterative feedback improvement. In such settings, the instance itself is only one part of the protocol. The textual description, environment interface, callable tools, feedback mechanism, and metrics must match the target capability being evaluated.

\subsection{Applicability of Evaluation Resources}

Different resources are suitable for different LLM roles, and dataset selection should match the capability under evaluation. 
For the modeler role, the emphasis is whether an LLM can identify decision variables, objective functions, and constraints from natural-language requirements,
 and then generate mathematical formulations or modeling code accepted by a solver. Resources containing business semantics, constraint descriptions, 
 and reference formulations are more appropriate for evaluating formalization ability.

For the designer role, the focus shifts to the effectiveness of generated heuristic rules, search operators, or executable solving programs. 
Standardized VRP benchmarks and LLM-oriented evaluation frameworks provide more unified instance generation, feasibility checking, and cost accounting, 
making them suitable for comparing solution quality and runtime efficiency after different models generate algorithms. For the coordinator role, 
the model often needs to organize external solvers, evaluators, repair modules, or neural solvers. 
Real or near-real operational datasets and interaction-based evaluation frameworks can better reflect tool invocation, 
workflow control, and multi-module collaboration.

\section{Comparative Experiments}
\label{sec:experiments}

Based on the datasets and evaluation frameworks discussed above, two complementary experiments are organized. 
The first uses CVRP as the basic test case and compares traditional solvers, NCO methods, and LLM-driven methods in terms of solution quality, 
solving time, and scalability. The second targets multiple VRP variants and evaluates the ability of different LLMs to generate VRP heuristics under 
the HeuriGym framework, focusing on cross-variant solution quality, valid returns under complex constraints, and the tradeoff between cost and solving time. 
All experiments were conducted in Python on a platform with an NVIDIA RTX 3090Ti GPU with 24 GB of memory and an Intel Core Ultra 7 265K CPU at 3.90 GHz.

In these experiments, the cost $C$ denotes the route objective value computed by the evaluator. For the VRP test data used here, 
$C$ mainly corresponds to the sum of travel distances over all vehicle routes. Capacity, time-window, backhaul-order, 
and route-length constraints are handled through feasibility checks. If a solution is infeasible, 
it is treated as invalid and is not included in the cost statistics.

\subsection{Solution-Quality Evaluation of Representative Methods}

\begin{table*}[!t]
\centering
\caption{Statistical results of CVRP experiments.}
\label{tab:cvrp_stat_results}
\footnotesize
\setlength{\tabcolsep}{2.6pt}
\renewcommand{\arraystretch}{1.12}
\begin{tabular*}{\textwidth}{@{\extracolsep{\fill}}lrrrrrrrrrr@{}}
\toprule
\multirow{2}{*}{Method} & \multicolumn{2}{c}{CVRP20} & \multicolumn{2}{c}{CVRP50} & \multicolumn{2}{c}{CVRP100} & \multicolumn{2}{c}{CVRP200} & \multicolumn{2}{c}{CVRP500}\\
\cmidrule(lr){2-3}\cmidrule(lr){4-5}\cmidrule(lr){6-7}\cmidrule(lr){8-9}\cmidrule(l){10-11}
& Gap\% ($\downarrow$) & Time (s) & Gap\% ($\downarrow$) & Time (s) & Gap\% ($\downarrow$) & Time (s) & Gap\% ($\downarrow$) & Time (s) & Gap\% ($\downarrow$) & Time (s)\\
\midrule
\rowcolor{tableTitleColor}
\multicolumn{11}{c}{Traditional solvers}\\
PyVRP(800) & 0.00\% & 33.89 & 0.00\% & 143.87 & 0.00\% & 345.44 & 0.00\% & 732.87 & 0.00\% & 2302.92\\
LKH3 & 0.17\% & 34.80 & 0.98\% & 573.29 & 1.92\% & 4344.31 & -- & -- & -- & --\\
OR-Tools (automatic) & 4.50\% & 1.26 & 8.41\% & 4.33 & 9.55\% & 14.66 & 9.91\% & 60.36 & 8.87\% & 421.02\\
\midrule
\rowcolor{tableTitleColor}
\multicolumn{11}{c}{NCO methods}\\
POMO(aug8x) & 28.40\% & $<1$ & 3.53\% & $<1$ & 0.45\% & $<1$ & 4.35\% & 2.97 & 20.25\% & 35.05\\
RouteFinder & 2.22\% & $<1$ & 1.15\% & $<1$ & 1.13\% & $<1$ & 3.26\% & 3.69 & 11.69\% & 46.81\\
Sym-NCO & 35.63\% & $<1$ & 2.74\% & $<1$ & 0.64\% & $<1$ & 3.61\% & 1.56 & 16.26\% & 9.09\\
\midrule
\rowcolor{tableTitleColor}
\multicolumn{11}{c}{LLM-driven methods}\\
From Words to Routes & 21.90\% & 28.20 & 30.92\% & 74.10 & 45.03\% & 114.62 & 55.64\% & 259.20 & -- & --\\
ReEvo & 13.87\% & 44.20 & 23.34\% & 96.62 & 31.32\% & 209.98 & 40.70\% & 502.99 & 43.78\% & 1938.65\\
G-LNS & 5.83\% & $<1$ & 8.37\% & 3.90 & 12.91\% & 25.96 & 16.78\% & 210.85 & 10.30\% & 3436.34\\
HeuriGym & 3.48\% & 2.35 & 3.68\% & 15.16 & 3.63\% & 4.14 & 5.12\% & 8.82 & 3.73\% & 569.64\\
HeurAgenix & 27.23\% & 6.02 & 27.35\% & 24.09 & 32.35\% & 226.67 & 34.97\% & 1600.47 & 29.54\% & 3692.68\\
NCO-LLMs (POMO-based) & 18.49\% & $<1$ & 4.02\% & $<1$ & 0.90\% & $<1$ & 1.04\% & 3.22 & 1.24\% & 36.25\\
\bottomrule
\end{tabular*}
\vspace{2pt}
\parbox{\textwidth}{\footnotesize Note: ``--'' indicates that no valid result was reported at this scale because the runtime exceeded 5000 s or no feasible solution was obtained after ten consecutive attempts.}
\end{table*}

\subsubsection{Experimental Settings}

The first experiment uses CVRP as the core evaluation task because its definition is clear, its constraint structure is fundamental, 
and mature benchmarks and reference solvers are available. Given a depot, customer coordinates, demands, and vehicle capacity, 
each method must output routes that start and end at the depot, visit every customer exactly once, and satisfy the capacity limit on each route. 
The test data consist of randomly generated CVRP instances at five scales: CVRP20, CVRP50, CVRP100, CVRP200, and CVRP500. 
These scales reveal how solution quality and runtime change as the number of customers increases.

The compared methods fall into three groups. Traditional baselines include PyVRP\footnote{\url{https://pyvrp.org/}}, LKH3\footnote{\url{http://webhotel4.ruc.dk/~keld/research/LKH-3/}}, and OR-Tools\footnote{\url{https://github.com/google/or-tools}}. PyVRP,
 built on high-performance heuristic frameworks such as HGS, is used as the reference baseline for gap computation. 
 LKH3 and OR-Tools represent a classical heuristic solver and an engineering-oriented open-source routing toolkit, respectively. 
 The NCO group includes POMO \cite{kwon2020pomo}, RouteFinder \cite{berto2025routefinder}, and Sym-NCO \cite{kim2022symnco}, 
 which test learned solvers on random in-distribution instances. The LLM-driven group includes From Words to Routes \cite{huang2024words}, 
 ReEvo \cite{ye2024reevo}, G-LNS \cite{zhao2025glns}, HeuriGym \cite{chen2025heurigym}, HeurAgenix \cite{yang2025heurageni}, 
 and NCO-LLMs \cite{tran2025llmnco}.

The main metrics are relative gap and solving time. The gap is defined as
\begin{equation}
  \mathrm{Gap}=\frac{C_{\mathrm{method}}-C_{\mathrm{ref}}}{C_{\mathrm{ref}}}\times 100\%
  \label{eq:gap}
\end{equation}
where $C_{\mathrm{method}}$ denotes the route cost generated by the evaluated method, and $C_{\mathrm{ref}}$ denotes the reference cost obtained by PyVRP on the same instance. All reported gaps are computed against the PyVRP reference. A negative gap is possible when an evaluated method finds a solution better than the PyVRP reference. Solving time only counts the algorithm execution stage and does not include LLM API calling time.

\subsubsection{Experimental Results}

Table~\ref{tab:cvrp_stat_results} summarizes solution quality and runtime across the five random CVRP scales. 
NCO methods and NCO-LLMs report GPU inference time, 
while traditional solvers and the other LLM-driven methods report CPU runtime.

\begin{table*}[!t]
\centering
\caption{Results of multi-variant heuristic generation experiments.}
\label{tab:multivariant_llm_results}
\footnotesize
\setlength{\tabcolsep}{2.2pt}
\renewcommand{\arraystretch}{1.12}
\begin{tabular*}{\textwidth}{@{\extracolsep{\fill}}lrrrrrrrrrrrr@{}}
\toprule
\multirow{2}{*}{Model} & \multicolumn{2}{c}{CVRP} & \multicolumn{2}{c}{SDVRP} & \multicolumn{2}{c}{VRPB} & \multicolumn{2}{c}{VRPTW} & \multicolumn{2}{c}{VRPL} & \multicolumn{2}{c}{OVRP}\\
\cmidrule(lr){2-3}\cmidrule(lr){4-5}\cmidrule(lr){6-7}\cmidrule(lr){8-9}\cmidrule(lr){10-11}\cmidrule(l){12-13}
& $C$ ($\downarrow$) & Time (s) & $C$ ($\downarrow$) & Time (s) & $C$ ($\downarrow$) & Time (s) & $C$ ($\downarrow$) & Time (s) & $C$ ($\downarrow$) & Time (s) & $C$ ($\downarrow$) & Time (s)\\
\midrule
\rowcolor{tableTitleColor}
\multicolumn{13}{c}{DeepSeek series}\\
deepseek-r1:14b & 850.56 & 2.26 & 921.55 & 2.27 & 997.15 & 2.25 & -- & -- & 894.84 & 2.22 & 589.73 & 2.22\\
deepseek-v4-flash & 693.76 & 3.19 & 775.59 & 2.17 & 805.42 & 2.11 & 1099.11 & 7.65 & 693.54 & 4.90 & 467.38 & 2.73\\
deepseek-v4-pro & 694.89 & 2.60 & 711.60 & 2.49 & 758.09 & 612.10 & 1064.67 & 611.12 & 693.86 & 4.78 & 434.97 & 351.50\\
\midrule
\rowcolor{tableTitleColor}
\multicolumn{13}{c}{Qwen series}\\
qwen3.5:9b & 1911.74 & 2.20 & 829.59 & 2.18 & -- & -- & 1589.82 & 2.19 & 894.84 & 2.17 & 625.04 & 6.34\\
qwen3.6:27b & 725.76 & 10.67 & 776.68 & 6.89 & 884.94 & 6.52 & 1075.92 & 7.04 & 781.13 & 3.17 & 480.17 & 9.42\\
qwen3.7-plus & 692.33 & 2.17 & 743.02 & 6.07 & 835.99 & 2.72 & 1080.62 & 22.26 & 695.90 & 2.18 & 432.84 & 2.44\\
\midrule
\rowcolor{tableTitleColor}
\multicolumn{13}{c}{GPT series}\\
gpt-oss:20b & 701.26 & 5.93 & 747.28 & 2.23 & 847.59 & 2.22 & 1205.83 & 5.02 & 695.71 & 4.64 & 443.80 & 6.37\\
gpt-5.4-mini & 769.67 & 2.00 & 803.52 & 1.96 & 856.36 & 2.24 & 1350.28 & 1.99 & 891.76 & 2.07 & 477.37 & 2.09\\
gpt-5.5 & 665.50 & 588.23 & 660.75 & 269.99 & 631.79 & 129.85 & 1032.09 & 619.10 & 667.11 & 607.17 & 419.85 & 553.76\\
\bottomrule
\end{tabular*}
\vspace{2pt}
\parbox{\textwidth}{\footnotesize Note: ``--'' indicates that the model did not obtain a valid solution after ten consecutive attempts.}
\end{table*}

The results show the expected contrast between traditional solvers and NCO methods. PyVRP, LKH3, and OR-Tools maintain strong solution quality, 
but runtime grows with the number of customers. LKH3 reports no valid result on CVRP200 and CVRP500 because of excessive runtime. 
POMO, RouteFinder, and Sym-NCO have short GPU inference times, yet their quality degrades on large instances: 
on CVRP500, their gaps reach 20.25\%, 11.69\%, and 16.26\%, respectively. Learned construction policies are fast at inference time, 
but their scalability remains tied to training distributions and model architectures.

LLM-driven methods behave differently depending on how the model is used. From Words to Routes directly generates routes from natural language, 
representing end-to-end route reasoning. Its gap increases from 21.90\% to 55.64\% as scale grows, 
and it fails to produce valid constrained routes on CVRP500. Direct route generation remains fragile in large combinatorial spaces where 
feasibility and quality must be maintained simultaneously. ReEvo generates executable solving programs, 
but its gap also rises from 13.87\% on CVRP20 to 43.78\% on CVRP500, 
indicating that automatically generated heuristics can lose effectiveness at larger scales.

External search, execution feedback, and candidate screening improve LLM-driven performance. G-LNS uses LLMs to generate search operators or neighborhood strategies and obtains gaps of 5.83\%--16.78\% from CVRP20 to CVRP500, generally outperforming direct route generation and generic automatic heuristic design. Its CVRP500 runtime, however, reaches 3436.34 s, showing the cost of search-enhanced improvement. HeuriGym keeps gaps within 3.48\%--5.12\% across all five scales, indicating that execution feedback, candidate selection, and feasibility checking stabilize heuristic quality. Runtime is not monotonic with instance size: CVRP50 takes longer than CVRP100 and CVRP200, which implies dependence on candidate-program quality, repair frequency, and early stopping. On CVRP500, HeuriGym still achieves a low gap of 3.73\%, but runtime rises to 569.64 s. The advantage of HeuriGym comes mainly from closed-loop quality control rather than the scale efficiency of a single generated heuristic. HeurAgenix remains relatively stable, with most gaps between 27\% and 35\%, but its overall accuracy is limited.

NCO-LLMs performs best on medium and large instances. Its gaps on CVRP100, CVRP200, and CVRP500 are 0.90\%, 1.04\%, and 1.24\%, lower than those of POMO, RouteFinder, and Sym-NCO at the same scales. It also takes only 36.25 s on CVRP500. These results support a hybrid pattern in which LLMs provide structural information extraction, strategy generation, or solver enhancement, while trained neural solvers or search frameworks carry out concrete solving.

\subsection{Evaluation of LLM-Generated Heuristics}

Beyond the single-task CVRP comparison, a multi-variant heuristic-generation experiment is designed to evaluate different LLMs under a unified workflow. The central question is whether an LLM can generate executable heuristic programs with acceptable solution quality when problem descriptions, instance data, and execution feedback are provided.

\subsubsection{Experimental Settings}

The multi-variant experiment uses RL4CO \cite{berto2025rl4co} to generate six VRP variants: CVRP, SDVRP, VRPB, VRPTW, VRPL, and OVRP. Each problem has 50 customers to control scale and preserve comparability. Each variant contains 16 training instances and 64 test instances, and reported results are computed on the test set. CVRP evaluates route construction under capacity constraints. SDVRP, VRPB, VRPTW, VRPL, and OVRP introduce split delivery, backhaul ordering, time windows, route-length limits, and open routes, respectively. These variants test whether LLMs can understand and implement variant-specific constraints.

The evaluated models include deepseek-r1:14b, deepseek-v4-flash, and deepseek-v4-pro from the DeepSeek series; qwen3.5:9b, qwen3.6:27b, and qwen3.7-plus from the Qwen series; and gpt-oss:20b, gpt-5.4-mini, and gpt-5.5 from the GPT series. The experiment uses the LLM-based heuristic-design framework HeuriGym \cite{chen2025heurigym}. HeuriGym is selected because the CVRP comparison shows relatively stable gaps across scales, and because its code generation, execution feedback, and candidate screening support comparison under the same problem descriptions, input-output formats, and evaluation process.

For each VRP variant, the prompt templates contain the problem constraints, input-output format, and feasibility requirements. All models generate solving code under the same prompt templates, instance sets, maximum iteration budget, and execution-time limit. A unified evaluator runs the generated programs, checks feasibility, and records route cost and solving time.

\subsubsection{Experimental Results}

\begin{figure*}[!t]
\centering
\includegraphics[width=0.82\textwidth]{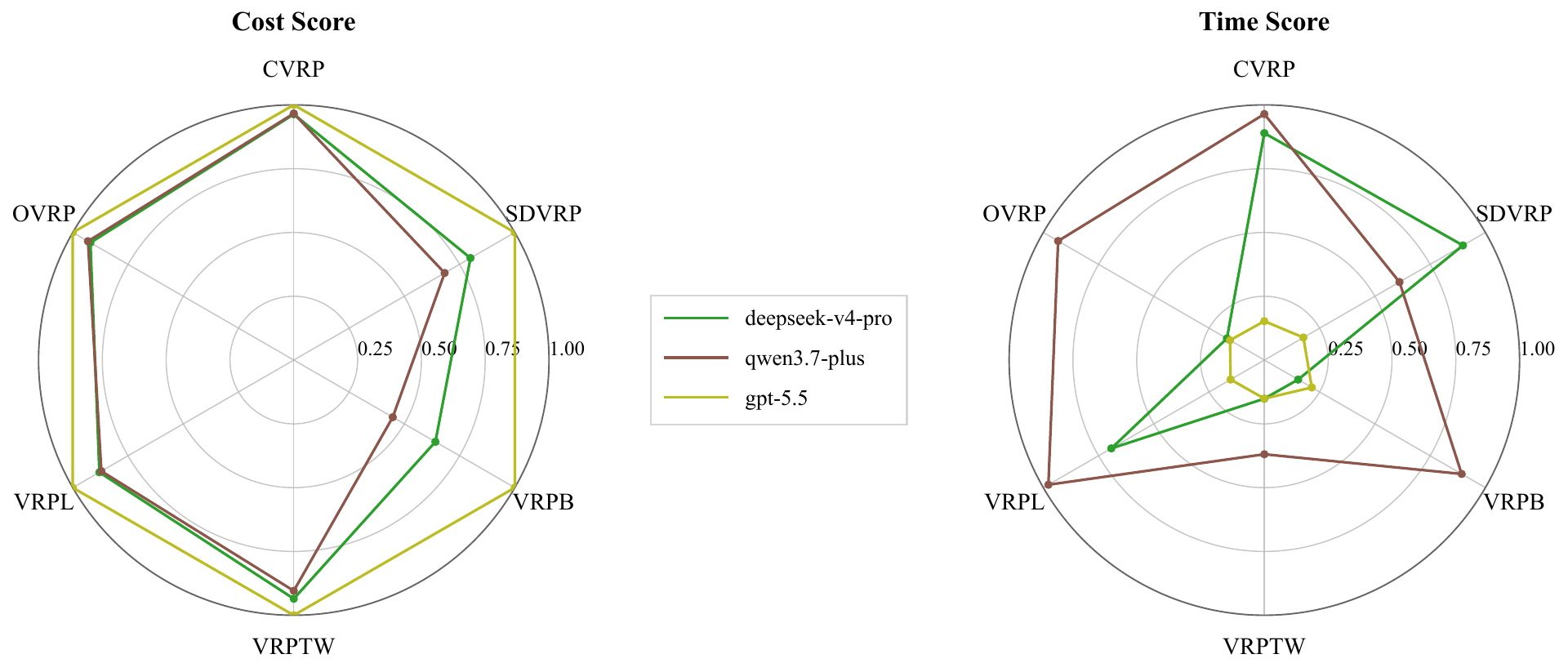}
\caption{Cost and time score comparison of LLM-generated heuristics across multiple VRP variants.}
\label{fig:llm_cost_log_radar}
\end{figure*}

Table~\ref{tab:multivariant_llm_results} reports heuristic-generation results across the six VRP variants. To compare solution quality and solving efficiency across variants, Fig.~\ref{fig:llm_cost_log_radar} plots normalized radar charts for solution cost and solving time. Larger normalized values indicate better performance. For solution cost, the transformation $\log(1+c_{m,p})$ is first applied to the cost $c_{m,p}$. Let
\begin{equation}
  \ell_{m,p}=\log(1+c_{m,p}), \quad
  D_p=\max_{m'}\ell_{m',p}-\min_{m'}\ell_{m',p}
  \label{eq:cost_aux}
\end{equation}
The reverse-normalized cost score within each problem variant $p$ is then computed as
\begin{equation}
  S^{\mathrm{cost}}_{m,p}=0.10+0.90\cdot
  \frac{\max_{m'}\ell_{m',p}-\ell_{m,p}}{D_p}
  \label{eq:cost_score}
\end{equation}
where $m$ denotes an LLM and $p$ denotes a VRP variant. This transformation assigns 1 to the lowest cost, 0.10 to the worst valid result, and 0 to missing or invalid results. For solving time, because most models finish within a few seconds while a few results take hundreds of seconds, square-root scaling relative to the fastest valid time is used:
\begin{equation}
  S^{\mathrm{time}}_{m,p}=0.10+0.90\cdot
  \left(\frac{\min_{m'}t_{m',p}}{t_{m,p}}\right)^{1/2}
  \label{eq:time_score}
\end{equation}
where $t_{m,p}$ is the solving time of model $m$ on variant $p$. This metric assigns 1 to the fastest valid result for each variant and 0 to missing or invalid results.

Across variants, constraint complexity first affects whether generated heuristics remain valid. CVRP and OVRP have relatively direct structures, and most models return feasible results. SDVRP and VRPL require the program to maintain remaining demand, vehicle capacity, and distance budget during construction. VRPB and VRPTW are harder: the former imposes linehaul-backhaul ordering, while the latter requires synchronized updates of arrival time, waiting time, and service windows. On these two strongly constrained tasks, smaller models such as deepseek-r1:14b and qwen3.5:9b show insufficient constraint understanding, incomplete state updates, or execution failure.

Among models that return results stably, solution quality and runtime show a clear tradeoff. gpt-5.5 obtains the lowest cost on all six variants: 665.50, 660.75, 631.79, 1032.09, 667.11, and 419.85 on CVRP, SDVRP, VRPB, VRPTW, VRPL, and OVRP, respectively. Its heuristics adapt well to different constraint forms, but runtime ranges from 129.85 s to 619.10 s, far above the 2--10 s range of most other models. deepseek-v4-pro shows a similar pattern: it improves over the flash version on SDVRP, VRPB, VRPTW, and OVRP, but takes 612.10 s, 611.12 s, and 351.50 s on VRPB, VRPTW, and OVRP. Lower cost is often obtained through more extensive search, repair, or feasibility maintenance.

Smaller LLMs tend to generate lightweight heuristics that return feasible solutions quickly. qwen3.7-plus obtains costs of 692.33, 743.02, 695.90, and 432.84 on CVRP, SDVRP, VRPL, and OVRP. deepseek-v4-flash returns valid results on all six tasks with short runtime, and gpt-oss:20b also solves most tasks at low computational cost. These fast strategies, however, lag behind lower-cost models on complex variants. On VRPB, gpt-5.5 achieves 631.79, while qwen3.7-plus, gpt-oss:20b, and gpt-5.4-mini obtain 835.99, 847.59, and 856.36. On VRPTW, gpt-5.4-mini takes only 1.99 s but reaches 1350.28, far worse than the slower gpt-5.5 and deepseek-v4-pro.

Complex constraints make differences among LLMs more visible in route-state updates, capacity handling, time-window reasoning, visit-order control, and search-strategy design. Larger models are more likely to generate low-cost heuristics, but they often rely on more extensive search or repair and take longer. Smaller models offer runtime advantages through lightweight heuristics, but solution quality can deteriorate on strongly constrained variants such as VRPB and VRPTW.

\section{Conclusion and Perspectives}
\label{sec:conclusion}

This paper surveyed LLM-driven VRP research through three roles: modeler, designer, and coordinator, corresponding to problem modeling, heuristic or route generation, and solving-workflow organization. It also reported CVRP benchmark comparisons and multi-variant heuristic-generation experiments to compare solution quality, runtime, scalability, and constraint handling across the main technical routes. Table~\ref{tab:llm_vrp_roles_summary} summarizes the main characteristics of the three roles.

\begin{table*}[!t]
\centering
\caption{Roles of LLMs in VRP solving.}
\label{tab:llm_vrp_roles_summary}
\footnotesize
\setlength{\tabcolsep}{5pt}
\renewcommand{\arraystretch}{1.16}
\begin{tabular*}{\textwidth}{@{\extracolsep{\fill}}>{\raggedright\arraybackslash}p{0.10\textwidth}>{\raggedright\arraybackslash}p{0.23\textwidth}>{\raggedright\arraybackslash}p{0.60\textwidth}@{}}
\toprule
Role & Function & Characteristics \\
\midrule
Modeler & Generate formulations, constraints, or modeling code &
Lowers the barrier to optimization modeling and supports rapid prototyping, but generated variables, indices, constraints, and objective functions must be checked carefully. \\
Designer & Generate heuristics, operators, or solving code &
Explores new construction and improvement strategies, while solution quality, feasibility, and code stability depend strongly on prompt design and execution feedback. \\
Coordinator & Organize solvers, evaluators, agents, and neural modules &
Improves tool collaboration and workflow automation, but requires clear system boundaries, controlled calling cost, error-propagation management, and reproducible protocols. \\
\bottomrule
\end{tabular*}
\end{table*}

LLMs provide an interface between natural-language business requirements and optimization pipelines. They can translate complex rules, generate heuristic ideas or executable programs, explain intermediate results, and invoke external tools, but they should complement rather than replace mature optimization methods. Their value is clearest when constraints are hard to formalize, variants change frequently, or fast solver prototyping is required. Future work should focus on the following directions.

\subsection{Interpretability and Verifiable Interfaces}

Constraints in practical VRP variants are often interdependent. Time windows affect waiting and service order, pickup-and-delivery rules modify vehicle states, and resource limits may interact with route duration or charging behavior. In such settings, final route cost alone cannot show whether an LLM has understood the constraints or why a generated heuristic fails. Future work should develop VRP-oriented intermediate representations and verifiable interfaces that convert natural-language tasks, mathematical models, heuristic code, and route outputs into checkable objects. Connecting these objects to solvers, simulators, or constraint checkers would improve both auditability and interpretability.

\subsection{Human-in-the-Loop Dispatching}

Real logistics dispatching often involves temporary orders, vehicle exceptions, traffic disruptions, and human operational experience. Plans may be revised repeatedly during execution. LLMs are well suited to serve as an interaction layer between dispatchers and optimization solvers by translating human intent into computable constraints, local re-optimization tasks, or solver-control instructions. Future systems should not only update routes, but also explain route changes and coordinate solution quality, business rules, and human preferences through multi-round feedback.

\subsection{Dataset Construction and Evaluation Protocols}

Standardized benchmarks and random instances are useful for reproducible experiments, but they cannot fully reflect road-network structures, traffic variation, operational rules, and human dispatching behavior. Future work should construct hierarchical evaluation resources in which standard instances, near-real road networks, and operational data serve different goals. Beyond solution quality and runtime, these resources should record process-level information such as code executability, constraint violations, repair traces, and failure reasons. Such information would make it possible to evaluate LLM generalization under changing rules and complex data distributions.

\subsection{Multimodal Routing Information}

Critical information in delivery tasks is not limited to text. Map structures, order records, historical trajectories, traffic states, and visual layouts can all affect routing decisions. Existing studies have begun to combine textual prompts with visual information for CVRP solving \cite{huang2024multimodalcvrp}, and to encode complex constraints through visual modalities before fusing them with graph representations for multitask VRP solving \cite{gui2026vafm}. Future research should examine how multimodal information contributes to constraint recognition, dynamic route adjustment, and solver control. Strict benchmarks and ablation studies are needed to verify whether multimodal systems provide real gains over text-only or graph-only approaches.

In summary, LLMs are expanding the design space of VRP research by connecting language understanding, code generation, search design, and tool orchestration. The next stage should move from feasibility evidence toward reliable, verifiable, and reproducible LLM-assisted routing systems. Progress will depend on closer coupling with mature optimization methods, clearer evaluation protocols, and datasets that expose both numerical performance and process-level behavior.

\bibliographystyle{IEEEtran}
\bibliography{IEEEabrv,references}

@article{dantzig1959truck,
  author  = {Dantzig, George B. and Ramser, John H.},
  title   = {The Truck Dispatching Problem},
  journal = {Management Science},
  volume  = {6},
  number  = {1},
  pages   = {80--91},
  year    = {1959},
  doi     = {10.1287/mnsc.6.1.80}
}

@article{laporte2009fifty,
  author={Laporte, Gilbert},
  title={Fifty Years of Vehicle Routing},
  journal={Transportation Science},
  volume={43},
  number={4},
  pages={408--416},
  year={2009},
  publisher={INFORMS},
  doi={10.1287/trsc.1090.0301}
}

@article{letchford2015mcf_cvrp,
  author={Letchford, Adam N. and Salazar-Gonz{\'a}lez, Juan Jos{\'e}},
  title={Stronger multi-commodity flow formulations of the {Capacitated Vehicle Routing Problem}},
  journal={European Journal of Operational Research},
  volume={244},
  number={3},
  pages={730--738},
  year={2015},
  publisher={Elsevier},
  doi={10.1016/j.ejor.2015.07.020}
}

@article{ropke2006adaptive,
  title={An Adaptive Large Neighborhood Search Heuristic for the Pickup and Delivery Problem with Time Windows},
  author={Ropke, Stefan and Pisinger, David},
  journal={Transportation Science},
  volume={40},
  number={4},
  pages={455--472},
  year={2006},
  publisher={INFORMS}
}

@article{solomon1987algorithms,
  author={Solomon, Marius M.},
  title={Algorithms for the Vehicle Routing and Scheduling Problems with Time Window Constraints},
  journal={Operations Research},
  volume={35},
  number={2},
  pages={254--265},
  year={1987},
  publisher={INFORMS}
}

@article{uchoa2017new,
  title={New Benchmark Instances for the Capacitated Vehicle Routing Problem},
  author={Uchoa, Eduardo and Pecin, Diego and Pessoa, Artur and Poggi, Marcus and Vidal, Thibaut and Subramanian, Anand},
  journal={European Journal of Operational Research},
  volume={257},
  number={3},
  pages={845--858},
  year={2017},
  publisher={Elsevier}
}

@article{archetti2006tabu,
  title={A Tabu Search Algorithm for the Split Delivery Vehicle Routing Problem},
  author={Archetti, Claudia and Hertz, Alain and Speranza, Maria Grazia},
  journal={Transportation Science},
  volume={40},
  number={1},
  pages={64--73},
  year={2006},
  publisher={INFORMS}
}

@article{schneider2014electric,
  title={The Electric Vehicle-Routing Problem with Time Windows and Recharging Stations},
  author={Schneider, Michael and Stenger, Andreas and Goeke, Dominik},
  journal={Transportation Science},
  volume={48},
  number={4},
  pages={500--520},
  year={2014},
  publisher={INFORMS},
  doi={10.1287/trsc.2013.0490}
}

@article{cordeau1997tabu,
  title={A Tabu Search Heuristic for Periodic and Multi-Depot Vehicle Routing Problems},
  author={Cordeau, Jean-Fran{\c{c}}ois and Gendreau, Michel and Laporte, Gilbert},
  journal={Networks},
  volume={30},
  number={2},
  pages={105--119},
  year={1997}
}

@article{goetschalckx1989vehicle,
  title={The Vehicle Routing Problem with Backhauls},
  author={Goetschalckx, Marc and Jacobs-Blecha, Catherine},
  journal={European Journal of Operational Research},
  volume={42},
  number={1},
  pages={39--51},
  year={1989},
  publisher={Elsevier}
}

@article{li2007open,
  title={The Open Vehicle Routing Problem: Algorithms, Large-Scale Test Problems, and Computational Results},
  author={Li, Feiyue and Golden, Bruce and Wasil, Edward},
  journal={Computers \& Operations Research},
  volume={34},
  number={10},
  pages={2918--2930},
  year={2007},
  publisher={Elsevier},
  doi={10.1016/j.cor.2005.11.018}
}

@incollection{christofides1979vehicle,
  title={The Vehicle Routing Problem},
  author={Christofides, Nicos and Mingozzi, Aristide and Toth, Paolo},
  booktitle={Combinatorial Optimization},
  editor={Christofides, Nicos and Mingozzi, Aristide and Toth, Paolo and Sandi, Claude},
  pages={315--338},
  year={1979},
  publisher={Wiley},
  address={Chichester}
}

@article{deepseek2025r1,
  title={{DeepSeek-R1}: Incentivizing Reasoning Capability in {LLMs} via Reinforcement Learning},
  author={{DeepSeek-AI}},
  journal={arXiv preprint arXiv:2501.12948},
  year={2025}
}

@article{hu2022lora,
  title={{LoRA}: Low-Rank Adaptation of Large Language Models},
  author={Hu, Edward J and Shen, Yelong and Wallis, Phillip and Allen-Zhu, Zeyuan and others},
  journal={arXiv preprint arXiv:2106.09685},
  year={2022}
}

@inproceedings{ouyang2022training,
  title={Training Language Models to Follow Instructions with Human Feedback},
  author={Ouyang, Long and Wu, Jeffrey and Jiang, Xu and Almeida, Diogo and Wainwright, Carroll and Mishkin, Pamela and Zhang, Chong and Agarwal, Sandhini and Slama, Katarina and Ray, Alex and others},
  booktitle={Advances in Neural Information Processing Systems (NeurIPS)},
  volume={35},
  pages={27730--27744},
  year={2022}
}

@article{dong2023sft,
  title={How Abilities in Large Language Models Are Affected by Supervised Fine-Tuning Data Composition},
  author={Dong, Guanting and Yuan, Hongyi and Lu, Kangjie and others},
  journal={arXiv preprint arXiv:2310.05492},
  year={2023}
}

@article{zhou2020progress,
  title={Progress in Neural {NLP}: Modeling, Learning, and Reasoning},
  author={Zhou, Ming and Duan, Nan and Liu, Shujie and Shum, Heung-Yeung},
  journal={Engineering},
  volume={6},
  number={3},
  pages={275--290},
  year={2020}
}

@article{min2023recent,
  title={Recent Advances in Natural Language Processing via Large Pre-Trained Language Models: A Survey},
  author={Min, Bonan and Ross, Hayley and Sulem, Elior and Veyseh, Amir Pouran Ben and Bhattacharjee, Abhik and Cattan, Arie and Yoon, Seunghyun and Kwon, Youngjae and Hsu, Yushi and others},
  journal={ACM Computing Surveys},
  volume={56},
  pages={1--40},
  year={2023},
  publisher={ACM}
}

@inproceedings{vinyals2015pointer,
  title={Pointer Networks},
  author={Vinyals, Oriol and Fortunato, Meire and Jaitly, Navdeep},
  booktitle={Advances in Neural Information Processing Systems (NeurIPS)},
  volume={28},
  year={2015}
}

@inproceedings{kool2018attention,
  title={Attention, Learn to Solve Routing Problems!},
  author={Kool, Wouter and van Hoof, Herke and Welling, Max},
  booktitle={International Conference on Learning Representations (ICLR)},
  year={2019}
}

@inproceedings{kwon2020pomo,
  title={{POMO}: Policy Optimization with Multiple Optima for Reinforcement Learning},
  author={Kwon, Yeong-Dae and Choo, Jinho and Kim, Byoungjip and Yoon, Iljoo and Gwon, Youngjune and Min, Seungjai},
  booktitle={Advances in Neural Information Processing Systems (NeurIPS)},
  volume={33},
  year={2020}
}

@inproceedings{kim2022symnco,
  title={{Sym-NCO}: Leveraging Symmetricity for Neural Combinatorial Optimization},
  author={Kim, Minsu and Park, Junyoung and Park, Jinkyoo},
  booktitle={Advances in Neural Information Processing Systems (NeurIPS)},
  volume={35},
  year={2022}
}

@inproceedings{luo2023lehd,
  title={Neural Combinatorial Optimization with Heavy Decoder: Toward Large Scale Generalization},
  author={Luo, Fu and Lin, Xi and Liu, Fei and Zhang, Qingfu and Wang, Zhenkun},
  booktitle={Advances in Neural Information Processing Systems (NeurIPS)},
  volume={36},
  year={2023}
}

@inproceedings{gasse2019gcnn,
  title={Exact Combinatorial Optimization with Graph Convolutional Neural Networks},
  author={Gasse, Maxime and Chetelat, Didier and Ferroni, Nicola and Charlin, Laurent and Lodi, Andrea},
  booktitle={Advances in Neural Information Processing Systems (NeurIPS)},
  volume={32},
  year={2019}
}

@inproceedings{hottung2020nlNS,
  title={Neural Large Neighborhood Search for the Capacitated Vehicle Routing Problem},
  author={Hottung, Andr{\'e} and Tierney, Kevin},
  booktitle={European Conference on Artificial Intelligence (ECAI)},
  pages={443--450},
  year={2020},
  doi={10.3233/FAIA200124}
}

@article{bengio2021machine,
  title={Machine Learning for Combinatorial Optimization: A Methodological Tour d'Horizon},
  author={Bengio, Yoshua and Lodi, Andrea and Prouvost, Antoine},
  journal={European Journal of Operational Research},
  volume={290},
  number={2},
  pages={405--421},
  year={2021},
  publisher={Elsevier}
}

@techreport{openai2023gpt4,
  title={{GPT-4} Technical Report},
  author={{OpenAI}},
  institution={OpenAI},
  year={2023}
}

@article{qwen2024qwen25,
  title={{Qwen2.5} Technical Report},
  author={Yang, An and Yang, Baosong and Zhang, Beichen and Hui, Binyuan and Zheng, Bo and Yu, Bowen and Li, Chengyuan and Liu, Dayiheng and Huang, Fei and Wei, Haoran and others},
  journal={arXiv preprint arXiv:2412.15115},
  year={2024}
}

@inproceedings{brown2020language,
  title={Language Models are Few-Shot Learners},
  author={Brown, Tom and Mann, Benjamin and Ryder, Nick and Subbiah, Melanie and Kaplan, Jared D and Dhariwal, Prafulla and Neelakantan, Arvind and others},
  booktitle={Advances in Neural Information Processing Systems (NeurIPS)},
  volume={33},
  pages={1877--1901},
  year={2020}
}

@article{chen2021evaluating,
  title={Evaluating Large Language Models Trained on Code},
  author={Chen, Mark and Tworek, Jerry and Jun, Heewoo and Yuan, Qiming and Pinto, Henrique Ponde de Oliveira and others},
  journal={arXiv preprint arXiv:2107.03374},
  year={2021}
}

@inproceedings{schick2023toolformer,
  title={Toolformer: Language Models Can Teach Themselves to Use Tools},
  author={Schick, Timo and Dwivedi-Yu, Jane and Dess{\`i}, Roberto and Raileanu, Roberta and Lomeli, Maria and others},
  booktitle={Advances in Neural Information Processing Systems (NeurIPS)},
  year={2023}
}

@inproceedings{vaswani2017attention,
  title={Attention Is All You Need},
  author={Vaswani, Ashish and Shazeer, Noam and Parmar, Niki and Uszkoreit, Jakob and Jones, Llion and Gomez, Aidan N and Kaiser, {\L}ukasz and Polosukhin, Illia},
  booktitle={Advances in Neural Information Processing Systems (NeurIPS)},
  volume={30},
  year={2017}
}

@inproceedings{liu2024eoh,
  title={Evolution of Heuristics: Towards Efficient Automatic Algorithm Design Using Large Language Model},
  author={Liu, Fei and Tong, Xialiang and Yuan, Mingxuan and Lin, Xi and Luo, Fu and Wang, Zhenkun and Lu, Zhichao and Zhang, Qingfu},
  booktitle={International Conference on Machine Learning (ICML)},
  year={2024}
}

@inproceedings{ye2024reevo,
  title={{ReEvo}: Large Language Models as Hyper-Heuristics with Reflective Evolution},
  author={Ye, Haoran and Wang, Jiarui and Cao, Zhiguang and Berto, Federico and Hua, Chuanbo and Kim, Haeyeon and Park, Jinkyoo and Song, Guojie},
  booktitle={Advances in Neural Information Processing Systems (NeurIPS)},
  year={2024}
}

@inproceedings{zheng2025mctsahd,
  title={Monte Carlo Tree Search for Comprehensive Exploration in {LLM}-Based Automatic Heuristic Design},
  author={Zheng, Zhi and Xie, Zhuoliang and Wang, Zhenkun and Hooi, Bryan},
  booktitle={International Conference on Machine Learning (ICML)},
  pages={78338--78373},
  year={2025}
}

@article{liu2024llm4ad,
  title={{LLM4AD}: A Platform for Algorithm Design with Large Language Model},
  author={Liu, Fei and Zhang, Rui and Xie, Zhuoliang and Tong, Xialiang and Yuan, Mingxuan and Zhang, Qingfu and Lin, Xi},
  journal={arXiv preprint arXiv:2412.17287},
  year={2024}
}

@inproceedings{ma2024autodh,
  title={Automatic Algorithm Design Assisted by {LLMs} for Solving Vehicle Routing Problems},
  author={Ma, Long and Hao, Xingxing and Yang, Ruikang and others},
  booktitle={International Conference on Signal Processing (ICSP)},
  pages={247--252},
  year={2024},
  doi={10.1109/ICSP62129.2024.10846261}
}

@inproceedings{dat2025hsevo,
  title={{HSEvo}: Elevating Automatic Heuristic Design with Diversity-Driven Harmony Search and Genetic Algorithm Using {LLMs}},
  author={Dat, Pham Van Tuan and Doan, Long and Binh, Huynh Thi Thanh},
  booktitle={AAAI Conference on Artificial Intelligence},
  pages={26931--26938},
  year={2025},
  doi={10.1609/aaai.v39i25.34898}
}

@inproceedings{yao2025meoh,
  title={Multi-Objective Evolution of Heuristic Using Large Language Model},
  author={Yao, Shengcai and Liu, Fei and Lin, Xi and others},
  booktitle={AAAI Conference on Artificial Intelligence},
  pages={27144--27152},
  year={2025},
  doi={10.1609/aaai.v39i25.34922}
}

@article{yang2025heurageni,
  title={{HeurAgenix}: Leveraging {LLMs} for Solving Complex Combinatorial Optimization Challenges},
  author={Yang, Xiaoxi and Zhang, Lei and Qian, Hui and others},
  journal={arXiv preprint arXiv:2506.15196},
  year={2025}
}

@article{huang2025calm,
  title={{CALM}: Co-evolution of Algorithms and Language Model for Automatic Heuristic Design},
  author={Huang, Ziyao and Wu, Weiwei and Wu, Kui and Wang, Jianping and Lee, Wei-Bin},
  journal={arXiv preprint arXiv:2505.12285},
  year={2025}
}

@article{li2025ars,
  title={{ARS}: Automatic Routing Solver with Large Language Models},
  author={Li, Kai and Liu, Fei and Wang, Zhuoyi and others},
  journal={arXiv preprint arXiv:2502.15359},
  year={2025}
}

@article{liu2023algorithm,
  title={Algorithm Evolution Using Large Language Model},
  author={Liu, Fei and Lin, Xi and Zhang, Qingfu and Tan, Kay Chen and Kwong, Sam},
  journal={arXiv preprint arXiv:2311.15249},
  year={2023}
}

@article{romeraparedes2024mathematical,
  title={Mathematical Discoveries from Program Search with Large Language Models},
  author={Romera-Paredes, Bernardino and Barekatain, Mohammadamin and Novikov, Alexander and Balog, Matej and Kumar, M Pawan and Dupont, Emilien and others},
  journal={Nature},
  volume={625},
  pages={468--475},
  year={2024},
  publisher={Nature Publishing Group}
}

@article{ahmaditeshnizi2023optimus,
  title={{OptiMUS}: Optimization Modeling Using {MIP} Solvers and Large Language Models},
  author={AhmadiTeshnizi, Ali and Gao, Wenzhi and Udell, Madeleine},
  journal={arXiv preprint arXiv:2310.06116},
  year={2023}
}

@inproceedings{zhang2024optllm,
  title={Solving General Natural-Language-Description Optimization Problems with Large Language Models},
  author={Zhang, Jihai and Wang, Wei and Guo, Siyan and Wang, Li and Lin, Fangquan and Yang, Cheng and Yin, Wotao},
  booktitle={the 2024 Conference of the North American Chapter of the Association for Computational Linguistics: Human Language Technologies (Volume 6: Industry Track)},
  pages={483--490},
  year={2024},
  doi={10.18653/v1/2024.naacl-industry.42}
}

@inproceedings{wang2025ormind,
  title={{ORM}ind: A Cognitive-Inspired End-to-End Reasoning Framework for Operations Research},
  author={Wang, Zhiyuan and Chen, Bokui and Huang, Yinya and Cao, Qingxing and He, Ming and Fan, Jianping and Liang, Xiaodan},
  booktitle={the 63rd Annual Meeting of the Association for Computational Linguistics (Volume 6: Industry Track)},
  pages={104--131},
  year={2025},
  doi={10.18653/v1/2025.acl-industry.10}
}

@article{huang2025orlm,
  title={{ORLM}: A Customizable Framework in Training Large Models for Automated Optimization Modeling},
  author={Huang, Chenyu and Tang, Zhengyang and Hu, Shixi and Jiang, Ruoqing and Zheng, Xin and Ge, Dongdong and Wang, Benyou and Wang, Zizhuo},
  journal={Operations Research},
  volume={73},
  number={6},
  pages={2986--3009},
  year={2025},
  doi={10.1287/opre.2024.1233}
}

@article{ma2024llamoco,
  title={{LLaMoCo}: Instruction Tuning of Large Language Models for Optimization Code Generation},
  author={Ma, Zeyuan and Guo, Hongshu and Chen, Jiacheng and Peng, Guojun and Cao, Zhiguang and Ma, Yining and Gong, Yue-Jiao},
  journal={arXiv preprint arXiv:2403.01131},
  year={2024}
}

@article{thind2025optimai,
  title={{OptimAI}: Optimizing Optimization Problems with {AI}},
  author={Thind, Rohan and Vashishtha, Vedant and Umenberger, Jack and Johansson, Karl Henrik},
  journal={arXiv preprint arXiv:2504.16918},
  year={2025}
}

@article{kong2025alphaopt,
  title={{AlphaOPT}: Formulating Optimization Programs with Self-Improving {LLM} Experience Library},
  author={Kong, Xin and Ning, Wenhao and Wu, Huayuan and Zheng, Yu},
  journal={arXiv preprint arXiv:2510.18428},
  year={2025}
}

@inproceedings{ding2026orr1,
  title={{OR-R1}: Automating Modeling and Solving of Operations Research Optimization Problem via Test-Time Reinforcement Learning},
  author={Ding, Zezhen and Tan, Zhen and Zhang, Jiheng and Chen, Tianlong},
  booktitle={the AAAI Conference on Artificial Intelligence},
  year={2026},
  doi={10.1609/aaai.v40i1.36983}
}

@inproceedings{jiang2025droc,
  title={{DRoC}: Elevating Large Language Models for Complex Vehicle Routing via Decomposed Retrieval of Constraints},
  author={Jiang, Xia and Wu, Yaoxin and Zhang, Chenhao and Zhang, Yingqian},
  booktitle={International Conference on Learning Representations (ICLR)},
  year={2025}
}

@article{yang2026orthought,
  title={{ORThought}: Benchmarking and Automating Logistics Optimization Modeling with Structured {LLM} Reasoning},
  author={Yang, Yichao and Li, Yuchen and Xiong, Zeyu and Tian, Haofei and Su, Qianwen and Wang, Xinyu and Wu, Dongxiao and He, Fei and Yin, Yafeng and Ye, Xiangyu and others},
  journal={Artificial Intelligence for Transportation},
  volume={6},
  pages={100059},
  year={2026},
  doi={10.1016/j.ait.2026.100059}
}

@article{huang2024words,
  title={From Words to Routes: Applying Large Language Models to Vehicle Routing},
  author={Huang, Zhehui and Shi, Guangyao and Sukhatme, Gaurav S.},
  journal={arXiv preprint arXiv:2403.10795},
  year={2024}
}

@inproceedings{ju2024ttg,
  title={{To the Globe (TTG)}: Towards Language-Driven Guaranteed Travel Planning},
  author={Ju, Da and Jiang, Song and Cohen, Andrew and Foss, Aaron and Mitts, Sasha and Zharmagambetov, Arman and Amos, Brandon and Li, Xian and Kao, Justine T. and Fazel-Zarandi, Maryam and Tian, Yuandong},
  booktitle={the 2024 Conference on Empirical Methods in Natural Language Processing: System Demonstrations},
  year={2024}
}

@inproceedings{ahmaditeshnizi2024optimus,
  title={{OptiMUS}: Scalable Optimization Modeling with {(MI)LP} Solvers and Large Language Models},
  author={Ahmaditeshnizi, Ali and Gao, Wenzhi and Udell, Madeleine},
  booktitle={International Conference on Machine Learning (ICML)},
  year={2024}
}

@inproceedings{liu2025optitree,
  title={{OptiTree}: Hierarchical Thoughts Generation with Tree Search for {LLM} Optimization Modeling},
  author={Liu, Haoyang and Wang, Jie and Cai, Yuyang and Han, Xiongwei and Kuang, Yufei and Hao, Jianye},
  booktitle={Advances in Neural Information Processing Systems (NeurIPS)},
  year={2025}
}

@article{zhou2025urs,
  title={{URS}: A Unified Neural Routing Solver for Cross-Problem Zero-Shot Generalization},
  author={Zhou, Changliang and Yu, Canhong and Yao, Shunyu and Lin, Xi and Wang, Zhenkun and Zhou, Yu and Zhang, Qingfu},
  journal={arXiv preprint arXiv:2509.23413},
  year={2025}
}

@article{liu2024llm_survey_algorithm_design,
  title={A Systematic Survey on Large Language Models for Algorithm Design},
  author={Liu, F. and Yao, Y. and Guo, P. and Yang, Z. and Zhao, Z. and Lin, X. and Tong, X. and Yuan, M. and Lu, Z. and Wang, Z. and others},
  journal={arXiv preprint arXiv:2410.14716},
  year={2024}
}

@article{xiu2026llmor_survey,
  title={Large Language Models for Operations Research: A Comprehensive Survey},
  author={Xiu, Xianchao and Li, Jianhao and Zhang, Jintao and Zheng, Yaoxin and Wang, Neng and Xue, Yu and Liu, Tiancheng and Tang, Jie and Li, Yong},
  journal={arXiv preprint arXiv:2605.20849},
  year={2026}
}

@article{wang2025llmor_survey,
  title={Large Language Models in Operations Research: Methods, Applications, and Challenges},
  author={Wang, Kai and Wan, Zhongwei and Jiao, Zhenyu and Chang, S. Y. and Ning, Jia and Zhang, Lingxin and Duan, Jin and Guo, Zhiqian and Zhou, Qingyu and Yang, Shiji and others},
  journal={arXiv preprint arXiv:2509.18180},
  year={2025}
}

@article{ghanbarzadeh2025structured,
  title={A Structured Review of Large Language Models in Metaheuristic Optimisation},
  author={Ghanbarzadeh, Shadi and Sumari, Putra and Vinuesa, Ricardo},
  journal={Machine Learning with Applications},
  volume={22},
  pages={100740},
  year={2025},
  doi={10.1016/j.mlwa.2025.100740}
}

@inproceedings{tran2025llmnco,
  title={Large Language Models powered Neural Solvers for Generalized Vehicle Routing Problems},
  author={Tran, Cong Dao and Nguyen-Tri, Quan and Binh, Huynh Thi Thanh and Thanh-Tung, Hoang},
  booktitle={ICLR 2025 Workshop on Towards Agentic AI for Science},
  year={2025}
}

@inproceedings{yao2022react,
  title={{ReAct}: Synergizing Reasoning and Acting in Language Models},
  author={Yao, Shunyu and Zhao, Jeffrey and Yu, Dian and Du, Nan and Shafran, Izhak and Narasimhan, Karthik and Cao, Yuan},
  booktitle={International Conference on Learning Representations (ICLR)},
  year={2023}
}

@misc{autogpt,
  title={{AutoGPT}},
  author={Richards, Toran Bruce},
  note={GitHub repository},
  year={2023}
}

@misc{langchain,
  title={{LangChain}: Building Applications with {LLMs} through Composability},
  author={Chase, Harrison},
  note={GitHub repository},
  year={2022}
}

@article{pisinger2007general,
  title={A general heuristic for vehicle routing problems},
  author={Pisinger, David and Ropke, Stefan},
  journal={Computers \& Operations Research},
  volume={34},
  number={8},
  pages={2403--2435},
  year={2007},
  publisher={Elsevier}
}

@article{vidal2013hybrid,
  title={A hybrid genetic algorithm with adaptive diversity management for a large class of vehicle routing problems with time windows},
  author={Vidal, Thibaut and Crainic, Teodor Gabriel and Gendreau, Michel and Prins, Christian},
  journal={Computers \& Operations Research},
  volume={40},
  number={1},
  pages={475--489},
  year={2013},
  publisher={Elsevier}
}

@inproceedings{zhang2025afl,
  title={An Agentic Framework with {LLMs} for Solving Complex Vehicle Routing Problems},
  author={Zhang, Ni and Cao, Zhiguang and Zhou, Jianan and Zhang, Cong and Ong, Yew-Soon},
  booktitle={International Conference on Learning Representations (ICLR)},
  year={2026}
}

@article{hottung2025vrpagent,
  title={{VRPAgent}: {LLM}-Driven Discovery of Heuristic Operators for Vehicle Routing Problems},
  author={Hottung, Andr{\'e} and Berto, Federico and Hua, Chuanbo and Ye, Haoran and Kim, Haeyeon and Park, Jinkyoo and Tierney, Kevin},
  journal={arXiv preprint arXiv:2510.07073},
  year={2025}
}

@inproceedings{zhu2025rfthgs,
  title={Refining Hybrid Genetic Search for {CVRP} via Reinforcement Learning-Finetuned {LLM}},
  author={Zhu, Rongjie and Zhang, Cong and Cao, Zhiguang},
  booktitle={the International Conference on Learning Representations (ICLR)},
  year={2026}
}

@article{berto2025routefinder,
  title={{RouteFinder}: Towards Foundation Models for Vehicle Routing Problems},
  author={Berto, Federico and Hua, Chuanbo and Zepeda, Nayeli Gast and Hottung, Andr{\'e} and Wouda, Niels and Lan, Leon and Park, Junyoung and Tierney, Kevin and Park, Jinkyoo},
  journal={Transactions on Machine Learning Research (TMLR)},
  year={2025}
}

@article{xie2025ailsahd,
  title={Enhancing {CVRP} Solver through {LLM}-Driven Automatic Heuristic Design},
  author={Xie, Zhuoliang and Liu, Fei and Wang, Zhenkun and Zhang, Qingfu},
  journal={arXiv preprint arXiv:2602.23092},
  year={2026}
}

@article{zhao2025glns,
  title={{G-LNS}: Generative Large Neighborhood Search for {LLM}-Based Automatic Heuristic Design},
  author={Zhao, Baoyun and Wang, He and Zeng, Liang},
  journal={arXiv preprint arXiv:2602.08253},
  year={2026}
}

@inproceedings{malik2026llmaide,
  title={{LLMAide}: Language-Assisted Neural Solver for Vehicle Routing Problems},
  author={Malik, Manuj and Zhou, Jianan and Jin, Yan and Cao, Zhiguang},
  booktitle={the International Conference on Autonomous Agents and Multiagent Systems (AAMAS) Extended Abstracts},
  year={2026},
  doi={10.65109/ISOA2063}
}

@article{zeng2026universalvrp,
  title={A Universal Framework for Vehicle Routing Problems with Large Language Models},
  author={Zeng, Rui-Bin and Yang, Mo-Xuan and Lei, Ming-Long and Niu, Ling-Feng and Dai, Yu-Hong},
  journal={Journal of the Operations Research Society of China},
  year={2026},
  doi={10.1007/s40305-026-00692-6}
}

@article{malik2026pyvrpplus,
  title={{PyVRP+}: {LLM}-Driven Metacognitive Heuristic Evolution for Hybrid Genetic Search in Vehicle Routing Problems},
  author={Malik, Manuj and Zhou, Jianan and Chirra, Shashank Reddy and Cao, Zhiguang},
  journal={arXiv preprint arXiv:2604.07872},
  year={2026}
}

@article{chi2026unslhe,
  title={A generalized neural solver based on {LLM}-guided heuristic evoluation framework for solving diverse variants of vehicle routing problems},
  author={Chi, Minyan and Pang, Wei and Wu, Xuan and Zhao, Peng and Li, Yuanshu and Wang, Tianfang and Qian, Junjie and Xiao, Yubin and Wang, Liupu and Zhou, You},
  journal={Expert Systems with Applications},
  volume={296},
  number={Part A},
  pages={128876},
  year={2026},
  doi={10.1016/j.eswa.2025.128876}
}

@inproceedings{chen2026dragon,
  title={{DRAGON}: {LLM}-Driven Decomposition and Reconstruction Agents for Large-Scale Combinatorial Optimization},
  author={Chen, Xinyun and Zhang, Xufang and Wu, Xing and Yang, Chengrun and Liao, Zhi and Li, Lijun and Xia, Shuo and Shi, Fusheng and Yang, Yuxiang and Ma, Wei},
  booktitle={the International Conference on Autonomous Agents and Multiagent Systems (AAMAS)},
  year={2026}
}

@article{feng2026alignopt,
  title={{AlignOPT}: Aligning Large Language Models with Graph Neural Solvers for Combinatorial Optimization},
  author={Feng, Shaodi and Lin, Zhuoyi and Wu, Yaoxin and Yin, Haiyan and Jin, Yan and Jayavelu, Senthilnath and Xu, Xun},
  journal={arXiv preprint arXiv:2603.27169},
  year={2026}
}

@article{wang2025maef,
  title={Multi-agent large language models as evolutionary optimizers for scheduling optimization},
  author={Wang, Yidan and Wang, Jiayin and Chu, Zhiwei},
  journal={Computers \& Industrial Engineering},
  volume={206},
  pages={111197},
  year={2025},
  doi={10.1016/j.cie.2025.111197}
}

@article{yan2026hmace,
  title={{HMACE}: Heterogeneous Multi-Agent Collaborative Evolution for Combinatorial Optimization},
  author={Yan, Dongming and Xue, Meiling and Wang, Zifan and Hao, Qingyao and Du, Wenbo},
  journal={arXiv preprint arXiv:2605.07214},
  year={2026}
}

@article{qu2026coral,
  title={{CORAL}: Towards Autonomous Multi-Agent Evolution for Open-Ended Discovery},
  author={Qu, Ao and Zheng, Han and Zhou, Zijian and Yan, Yihao and Tang, Yihong and Ong, Shao Yong and Hong, Fenglu and Zhou, Kaichen and Jiang, Chonghe and Kong, Minwei and Zhu, Jiacheng and Jiang, Xuan and Li, Sirui and Wu, Cathy and Low, Bryan Kian Hsiang and Zhao, Jinhua and Liang, Paul Pu},
  journal={arXiv preprint arXiv:2604.01658},
  year={2026}
}

@inproceedings{cao2025llmqlearningcvrptw,
  title={A Large Language Model-Enhanced {Q-Learning} for Capacitated Vehicle Routing Problem with Time Windows},
  author={Cao, Linjiang and Wang, Maonan and Xiong, Xi},
  booktitle={2025 IEEE 28th International Conference on Intelligent Transportation Systems (ITSC)},
  year={2025}
}

@article{shi2026llmvd,
  title={{LLM}-based automatic heuristic design for vehicle-drone collaborative routing problems},
  author={Shi, Haiyang and Zhen, Lu},
  journal={Transportation Research Part E: Logistics and Transportation Review},
  volume={209},
  pages={104760},
  year={2026},
  doi={10.1016/j.tre.2026.104760}
}

@inproceedings{shi2026moh,
  title={Generalizable Heuristic Generation Through {LLMs} with Meta-Optimization},
  author={Shi, Yiding and Zhou, Jianan and Song, Wen and Bi, Jieyi and Wu, Yaoxin and Cao, Zhiguang and Zhang, Jie},
  booktitle={the International Conference on Learning Representations (ICLR)},
  year={2026}
}

@inproceedings{zhu2026evoreal,
  title={Bridging Synthetic and Real Routing Problems via {LLM}-Guided Instance Generation and Progressive Adaptation},
  author={Zhu, Jianghan and Wu, Yaoxin and Lin, Zhuoyi and Zhang, Zhengyuan and Yin, Haiyan and Cao, Zhiguang and Jayavelu, Senthilnath and Li, Xiaoli},
  booktitle={the AAAI Conference on Artificial Intelligence},
  volume={40},
  pages={36591--36599},
  year={2026}
}

@inproceedings{zafar2025evrp,
  title={Using Large Language Models to Solve the Electric Vehicle Routing Problem with Advanced Prompting Techniques},
  author={Zafar, Usman and Bayhan, Sertac},
  booktitle={2025 IEEE 19th International Conference on Compatibility, Power Electronics and Power Engineering (CPE-POWERENG)},
  year={2025},
  doi={10.1109/CPE-POWERENG63314.2025.11027311}
}

@article{merchan2022amazonlastmile,
  title={2021 {Amazon} Last Mile Routing Research Challenge: Data Set},
  author={Merch{\'a}n, Daniel and Arora, Jatin and Pachon, Julian and Konduri, Karthik and Winkenbach, Matthias and Parks, Steven and Noszek, Joseph},
  journal={Transportation Science},
  year={2022},
  doi={10.1287/trsc.2022.1173}
}

@article{vrani2025deliveringdata,
  title={Delivering Data: A Real-World Dataset for Last-Mile Delivery Optimization},
  author={Vrani, Anna and Apostolidis, Savvas D. and Kapoutsis, Athanasios Ch. and Kosmatopoulos, Elias B.},
  journal={Data in Brief},
  volume={61},
  pages={111762},
  year={2025},
  doi={10.1016/j.dib.2025.111762}
}

@article{greenberg2025earli,
  title={Accelerating Vehicle Routing via {AI}-Initialized Genetic Algorithms},
  author={Greenberg, Ido and Sielski, Piotr and Linsenmaier, Hugo and Gandham, Rajesh and Mannor, Shie and Fender, Alex and Chechik, Gal and Meirom, Eli},
  journal={arXiv preprint arXiv:2504.06126},
  year={2025}
}

@article{heakl2025svrpbench,
  title={{SVRPBench}: A Realistic Benchmark for Stochastic Vehicle Routing Problem},
  author={Heakl, Ahmed and Shaaban, Yahia Salaheldin and Lahlou, Salem and Tak{\'a}{\v{c}}, Martin and Iklassov, Zangir},
  journal={arXiv preprint arXiv:2505.21887},
  year={2025}
}

@inproceedings{berto2025rl4co,
  title={{RL4CO}: an Extensive Reinforcement Learning for Combinatorial Optimization Benchmark},
  author={Berto, Federico and Hua, Chuanbo and Park, Junyoung and Luttmann, Laurin and Ma, Yining and Bu, Fanchen and Wang, Jiarui and Ye, Haoran and Kim, Minsu and Choi, Sanghyeok and Zepeda, Nayeli Gast and Hottung, Andr{\'e} and Zhou, Jianan and Bi, Jieyi and Hu, Yu and Liu, Fei and Kim, Hyeonah and Son, Jiwoo and Kim, Haeyeon and Angioni, Davide and Kool, Wouter and Cao, Zhiguang and Zhang, Jie and Shin, Kijung and Wu, Cathy and Ahn, Sungsoo and Song, Guojie and Kwon, Changhyun and Xie, Lin and Park, Jinkyoo},
  booktitle={the 31st {ACM} {SIGKDD} Conference on Knowledge Discovery and Data Mining},
  year={2025}
}

@article{chen2025heurigym,
  title={{HeuriGym}: An Agentic Benchmark for {LLM}-Crafted Heuristics in Combinatorial Optimization},
  author={Chen, Hongzheng and Wang, Yingheng and Cai, Yaohui and Hu, Hins and Li, Jiajie and Huang, Shirley and Deng, Chenhui and Liang, Rongjian and Kong, Shufeng and Ren, Haoxing and Samaranayake, Samitha and Gomes, Carla P. and Zhang, Zhiru},
  journal={arXiv preprint arXiv:2506.07972},
  year={2025}
}

@article{sun2025cobench,
  title={{CO-Bench}: Benchmarking Language Model Agents in Algorithm Search for Combinatorial Optimization},
  author={Sun, Weiwei and Feng, Shengyu and Li, Shanda and Yang, Yiming},
  journal={arXiv preprint arXiv:2504.04310},
  year={2025}
}

@article{albalkhi2026routeoptimization,
  title={Route Optimization Reimagined: Multi-Modal Large Language Models for Next-Generation Vehicle Routing},
  author={Albalkhi, Sireen Y. and Alotaibi, Deem F. and Dimitriou, Tassos and Ahmad, Imtiaz},
  journal={IEEE Access},
  year={2026},
  doi={10.1109/ACCESS.2026.3663141}
}

@article{huang2024multimodalcvrp,
  title={How Multimodal Integration Boost the Performance of {LLM} for Optimization: Case Study on Capacitated Vehicle Routing Problems},
  author={Huang, Yuxiao and Zhang, Wenjie and Feng, Liang and Wu, Xingyu and Tan, Kay Chen},
  journal={arXiv preprint arXiv:2403.01757},
  year={2024}
}

@article{gui2026vafm,
  title={Vision-Assisted Foundation Model for Solving Multi-Task Vehicle Routing Problems},
  author={Gui, Shuangchun and Cao, Zhiguang and Song, Wen and Ong, Yew-Soon},
  journal={arXiv preprint arXiv:2606.10431},
  year={2026}
}

\end{document}